\documentclass{amsart}[12pt]
\usepackage{amsthm}
\usepackage{amsfonts}
\usepackage{comment}
\usepackage{amssymb}
\usepackage{amsmath}
\usepackage{xcolor}
\usepackage{eufrak}
\usepackage{comment}
\usepackage[colorlinks = true, citecolor = blue]{hyperref}
\usepackage{mathtools}
\usepackage{eucal}
\usepackage{tikz-cd}
\numberwithin{equation}{section}
\numberwithin{figure}{section}
\usepackage{breqn}
\allowdisplaybreaks
\usepackage[OT2,T1]{fontenc}
\usepackage{cleveref}
\usepackage[shortlabels]{enumitem}
\usepackage{stmaryrd}
\newtheorem{theorem}{Theorem}

\newtheorem{proposition}[theorem]{Proposition}

\newtheorem{Theorem}{Theorem}[section]
\newtheorem{Corollary}[Theorem]{Corollary}
\newtheorem{Lemma}[Theorem]{Lemma}
\newtheorem{Proposition}[Theorem]{Proposition}
\theoremstyle{definition}
\newtheorem{Example}[Theorem]{Example}
\newtheorem{Definition}[Theorem]{Definition}
\theoremstyle{remark}
\newtheorem{remark}[Theorem]{Remark}
\newtheorem*{question}{Question}
\usepackage[margin=1.25in]{geometry}
\usepackage{cite}
\usepackage{tikz}
\usepackage{float}

\DeclareMathOperator{\Gal}{Gal}

\DeclareMathOperator{\Disc}{Disc}
\DeclareMathOperator{\Frob}{Frob}

\DeclareMathOperator{\lcm}{lcm}

\DeclareMathOperator{\NP}{NP}
\DeclareMathOperator{\Pic}{Pic}

\DeclareMathOperator{\rk}{rk}
\DeclareMathOperator{\height}{ht}
\DeclareMathOperator{\Sym}{Sym}

\newcommand{\Q}{\mathbb{Q}}

\newcommand{\C}{\mathbb{C}}
\newcommand{\Z}{\mathbb{Z}}
\newcommand{\R}{\mathbb{R}}

\renewcommand{\P}{\mathbb{P}}
\renewcommand{\a}{\boldsymbol{a}}
\newcommand{\balpha}{\boldsymbol{\alpha}}
\renewcommand{\b}{\boldsymbol{b}}

\newcommand{\wt}[1]{\widetilde{#1}}
\newcommand{\set}[1]{\left\lbrace #1 \right\rbrace}
\newcommand{\irr}{\mathrm{irr}}

\title{Fields generated by points on superelliptic curves}

\author{Lea Beneish}
\address{Lea Beneish\\
Department of Mathematics\\
University of North Texas\\
Denton, TX\\
United States}
\email{lea.beneish@unt.edu}

\author{Christopher Keyes} 
\address{Christopher Keyes, Department of Mathematics, King's College London, London, United Kingdom and Heilbronn Institute for Mathematical Research, Bristol, United Kingdom}
\email{christopher.keyes@kcl.ac.uk}

\makeatletter
\@namedef{subjclassname@1991}{\emph{2020} Mathematics Subject Classification}
\makeatother


\begin{document}

\begin{abstract} 
We give an asymptotic lower bound on the number of field extensions generated by algebraic points on superelliptic curves over $\Q$ with fixed degree $n$ and discriminant bounded by $X$. For $C$ a fixed such curve given by an affine equation $y^m = f(x)$ where $m \geq 2$ and $d= \deg f (x) \geq m$, we find that for all degrees $n$ divisible by $\gcd(m, d)$ and sufficiently large, the number of such fields is asymptotically bounded below by $X^{\delta_n}$, where $\delta_n \to 1/m^2$ as $n \to \infty$. We then give geometric heuristics suggesting that for n not divisible by $\gcd(m, d)$, degree $n$ points may be less abundant than those for which $n$ is divisible by $\gcd(m,d)$ and provide an example of conditions under which a curve is known to have finitely many points of certain degrees.
\end{abstract}

\maketitle

\section{Introduction}
Let $K$ be a number field, and let $C/K$ be a smooth curve of genus $g$. Faltings
\cite{faltings} proved that when $g\geq 2$, the set of $K$-rational points on $C$, $C(K)$, is finite, and in fact $C(L)$ is finite for any finite extension $L/K$. It is natural to ask if similar finiteness results hold for the higher degree points of $C$. We say the degree of an algebraic point $P \in C(\overline{K})$ is the degree $[K(P):K]$, where $K(P)$ is the minimal field of definition for $P$. While in fact a curve of genus $g \geq 2$ may have infinitely many points of some degree $n > 1$, it is still an interesting problem to characterize when this occurs and to prove finiteness results for ``sporadic'' points.
There have been several recent works related to the study of higher degree points on families of hyperelliptic curves (see \cite{BGW, GM}) and on various modular curves (see \cite{BourdonVirayEtAl, Box, BGRW, BN, DEvHMZB, OS}).

Instead of studying the points of $C$, one can take the perspective of studying the set of field extensions $K(P)/K$ generated by algebraic points $P \in C(\overline{K})$. This idea was suggested by Mazur and Rubin \cite{mazurrubin} in their Diophantine stability program, where a variety over $K$ is said to be Diophantine stable for $L/K$ if its $K$-rational points and $L$-rational points coincide. A natural first question to ask is: how many extensions generated by an algebraic point exist for a fixed degree when ordered by discriminant?

 Fixing the base field $K=\Q$, we define the following functions for counting  number fields by discriminant. Let 
\[N_n(X) = \#\left\lbrace L/\Q : [L:\Q] = n, \ \left|\Disc L/\Q\right| \leq X \right\rbrace ,\]
where $X > 0$ is a real number and $n \geq 1$ is any positive integer. For a fixed curve $C/\Q$, we define the counting function for extensions generated by an algebraic point of $C$ to be
\[N_{n, C}(X) = \#\left\lbrace \Q(P)/\Q : P \in C(\overline{\Q}), \ [\Q(P):\Q] = n, \ \left|\Disc \Q(P)/\Q\right| \leq X \right\rbrace.\]
 We further define 
\[N_{n, C}(X,G) = \#\left\lbrace \Q(P)/\Q : P \in C(\overline{\Q}), \ [\Q(P):\Q] = n, \ \left|\Disc \Q(P)/\Q\right| \leq X,  \ \Gal(\wt{\Q(}P)/\Q) \simeq G \right\rbrace\]
where $G$ is a permutation subgroup of the symmetric group $S_n$ and $\wt{\Q(}P)$ denotes the Galois closure of $\Q(P)/\Q$. 

When $E$ is an elliptic curve over $\Q$, Lemke Oliver and Thorne \cite{LOT} show $N_{n,E/\Q}(X, S_n) \gg X^{c_n-\epsilon}$ for a positive constant $c_n$ approaching 1/4 from below as $n \to \infty$. Conditionally, this exponent can be improved to approach $1/4$ from above. In fact, they show something stronger, namely that $X^{c_n - \epsilon}$ is an asymptotic lower bound on degree $n$ extensions for which the Mordell--Weil ranks satisfy $\rk E(K) > \rk E(\Q)$, with  specified root number. In \cite{keyes} the second author proves that for a hyperelliptic curve $C/\Q$ of genus $g\geq 1$ and for $n$ sufficiently large relative to $C$, with $n$ even if the defining polynomial of $C$ has even degree, we have $N_{n,C/\Q}(X,S_n) \gg X^{c_n}$, where $c_n$ is again a constant depending on $g$ which tends to $1/4$ from below as $n \to \infty$. 

We continue this program of studying the the set of fields generated by points on curves defined over $\Q$ in the case of superelliptic curves. For a positive integer $m \geq 2$, a \textit{superelliptic curve $C/\Q$} is a smooth projective geometrically integral curve given by an affine equation of the form
\begin{equation}\label{eq:sec}
	C \colon y^m = f(x) = \sum_{i=0}^d c_ix^i,
\end{equation}
where $f(x) \in \Z[x]$ is a polynomial of degree $d$ such that 
\begin{equation}\label{eq:f_not_power}
    f \text{ is } m \text{-th power free and } f \notin \left(\overline{\Q}[x]\right)^{m'} \text{ for all } m' \mid m,\ m' \neq 1.
\end{equation}
The condition \eqref{eq:f_not_power} on $f$ is satisfied if and only if the curve $C$ is geometrically integral. In this paper, we restrict further to the case where $m \leq d$. Such a curve $C$ possesses a degree $m$ map to the projective line $\P^1$ defined over $\Q$, sending a point $(x,y) \mapsto x$.  When $\gcd(m,d) \mid n$ and $n$ is sufficiently large, we have the following asymptotic lower bound for $N_{n,C}(X)$. 

\begin{theorem}\label{thm:bound}
Fix integers $m \geq 2$, $d \geq m$, a polynomial $f \in \Z[x]$ of degree $d$ satisfying \eqref{eq:f_not_power}, and let $C$ be the curve with equation \eqref{eq:sec}. If $n_0 = \max\left(d, \lcm(m,d) - m - d + 1\right)$, then for all $n \geq n_0$ such that $\gcd(m,d) \mid n$, subject to the additional constraint that $n \geq \frac{m^2(m-1)}{2}$ when $1 < \gcd(m,n) < m$, we have
\begin{equation}\label{eq:final_bound}
	N_{n,C}(X) \gg X^{\delta_n},
\end{equation}
for a constant $\delta_n$ depending on $m, d,$ and $n$ given explicitly in \eqref{eq:c_n}. The implied constant in \eqref{eq:final_bound} depends on $n$ and (the equation for) $C$.

Moreover, for all sufficiently large $n$ (relative to $m$ and $d$) with $\gcd(m,d) \mid n$, we have $N_{n,C}(X) \gg X^{\delta_n'}$ where
\begin{equation}\label{eq:final_improved}
	\delta_n' ={\frac{1}{m^2} \left(1 + \frac{(2m - 2dr + 1)n + d^2r^2 - mdr + mk - k^2}{2n(n-1)}\right)},
\end{equation}
with $1 \leq r < m$ and $0 \leq k < m$ are integers depending only on the residue classes of $n,d \pmod{m}$.
\end{theorem}

\begin{remark}
	We make note of a few properties of the constant $\delta_n$ in Theorem \ref{thm:bound}.
	\begin{enumerate}[label = (\roman*)]
		\item For any fixed choice of $m, d$, the constant $\delta_n$ in \eqref{eq:final_bound} satisfies $\delta_n - \frac{1}{m^2} \sim \frac{m - m^2 - dr + 3}{m^2(n-1)}$ in the limit as $n \to \infty$, where $1 \leq r < m$ is an integer depending only on $n,d \pmod{m}$. In particular, $\frac{m - m^2 - dr + 3}{m^2(n-1)}$ is negative, so we can say that in \eqref{eq:final_bound}, $\delta_n$ approaches $\frac{1}{m^2}$ from below.
		
		\item In contrast, the improved exponent in \eqref{eq:final_improved} satisfies $\delta_n' - \frac{1}{m^2} \sim \frac{2m - 2dr + 1}{2m^2(n-1)}$. In the case $m=d$ we have $r=1$ and thus $2m - 2dr + 1 = 1$, so $\delta_n \to \frac{1}{m^2}$ from above as $n \to \infty$. If $m < d$, the improved $\delta_n'$ will approach $\frac{1}{m^2}$ from below as in \eqref{eq:final_bound}.
		
		\item The improved exponent in \eqref{eq:final_improved} takes effect when we have good enough asymptotic upper bounds for $N_n(X)$. The bound of Lemke Oliver--Thorne \cite[Theorem 1.1]{LOT_upper} suffices when $n$ is taken to be large. We discuss for which $n$ the bound $N_{n,C}(X) \gg X^{\delta_n'}$ is known to hold in Section \ref{sec:improve}; see Figure \ref{fig:improved_bound}.
		
		\item Theorem \ref{thm:bound} agrees with or improves upon known lower bounds for $N_{n,C}(X)$ in the cases where $C$ is an elliptic curve \cite{LOT} or a hyperelliptic curve \cite{keyes}, though notably it does not specify the Galois group. This is partially addressed for the $m=3$ case in Theorem \ref{thm:m=3_bound}.
		
		\item We do not expect this lower bound to be sharp; in the case where $C$ is an elliptic curve, Lemke Oliver--Thorne \cite{LOT} suggest a heuristic of $X^{3/4 + o(1)}$ for the asymptotics of the number of fields $K/\Q$ for which $\rk E(K) = \rk E(\Q)+2$.
	\end{enumerate}
\end{remark}

The strategy for proving Theorem \ref{thm:bound}, employed also in \cite{LOT} and \cite{keyes}, is to use the equation for $C/\Q$ to find an explicit parameterized family of polynomials generating degree $n$ extensions $\Q(P)/\Q$. Some effort is required to verify that the members of the family are in fact irreducible. We then count the polynomials in this family and bound how often the number fields they generate are isomorphic.

It is more difficult to generalize the approach when it comes to showing that the polynomials generated have Galois group $S_n$. However, we offer the following partial result for counting fields with largest possible Galois group generated by points on trigonal superelliptic curves.

\begin{theorem}\label{thm:m=3_bound}
    Suppose $m=3$ and $3 \nmid d$. Let $C \colon y^3 =f(x)$ for $f$ a degree $d$ polynomial satisfying \eqref{eq:f_not_power}. Then for all integers $n \geq \max(2d - 2, 14)$ such that $n \equiv 2,4 \pmod{6}$, we have
    \begin{align*}
        N_{n,C}(X,S_n) &\gg X^{\delta_n} \text{ and}\\ 
        N_{n,C}(X,S_n) &\gg X^{\delta_n'} \text{ when } n \text{ is sufficiently large},
    \end{align*}
    with $\delta_n, \delta_n'$ as in Theorem \ref{thm:bound}.
\end{theorem}

Another notable limitation of Theorem \ref{thm:bound} is the condition that the count only applies for field extensions of degree $n$ where $n$ is such that $\gcd(m, d) \mid n$. In the case where $C$ is a hyperelliptic curve, we have $m=2$ and $d$ can be chosen to be odd if and only if $C$ has a rational Weierstrass point. In this case, $\gcd(m,d) = 1$, and our parameterization produces infinite families of odd degree $n$ points for $n$ sufficiently large. In the general case however, we take $d = 2g + 2$, where $g$ is the genus of $C$, giving $\gcd(m, d)=2$, so this parametrization does not produce any odd degree points (cf. \cite{keyes}). This is consistent with a result of Bhargava--Gross--Wang \cite{BGW} that says a positive proportion of locally soluble hyperelliptic curves have no odd degree points. In Section \ref{sec:failsafe} we speculate as to whether for superelliptic curves, points of degrees $n$ such that $\gcd(m, d) \mid n$ are ``more common” than points of degrees $n$ where $\gcd(m, d) \nmid n$. This section contains a description of various geometric sources from which we expect to find infinitely many points on these curves. We also discuss the relationship of these sources to the points obtained by the parameterization strategy. As a first step towards making these heuristics concrete, we prove the following proposition.
\begin{proposition}\label{prop:degN_2adic}
		Suppose $m, d$ are positive even integers such that $d > 4$. Let $N < \frac{d}{2} - 1$ have $2$-adic valuation strictly less than that of $m$, i.e.\ $v_2(N) < v_2(m)$. Then for a positive proportion approaching 100\% of squarefree degree $d$ polynomials $f(x)$, ordered by height, the superelliptic curve $C \colon y^m = f(x)$ has finitely many points of degree $N$.
\end{proposition}

This paper is organized as follows. In Section \ref{sec:parametrization} we give an overview of the parameterization strategy used in the proof of Theorems \ref{thm:bound} and \ref{thm:m=3_bound}. In Section \ref{sec:galois_newton} we recall the Newton polygon of a polynomial, and how it may be used to identify cycles in its Galois group over $\Q_p$. We then recall a criterion for a transitive permutation group to be the full symmetric group, based on containing cycles of certain lengths. Section \ref{sec:irred_sn} is devoted to proving that our parameterization strategy almost always produces irreducible polynomials, then in Section \ref{sec:m=3} we specialize to $m=3$ and prove that these polynomials have Galois group $S_n$ in certain cases. In Section \ref{sec:lower_bounds} we count polynomials produced by our parameterization and adjust for multiplicity to obtain a lower bound for $N_{n,C}(X)$, completing the proofs of Theorems \ref{thm:bound} and \ref{thm:m=3_bound}. A discussion of the geometric sources for infinite collections of points on superelliptic curves, and their relevance to field counting problems of this flavor, is given in Section \ref{sec:failsafe}.

\section*{Acknowledgments}
The authors are grateful to Henri Darmon, Hannah Larson, Robert Lemke Oliver, Dino Lorenzini, Jackson Morrow, Frank Thorne, Brooke Ullery, Isabel Vogt, and David Zureick-Brown for helpful conversations. The authors would further like to thank Abbey Bourdon, Hannah Larson, Robert Lemke Oliver, Dino Lorenzini, Jackson Morrow, Jeremy Rouse, Isabel Vogt, and the anonymous referee for their thoughtful comments on an earlier draft. 

CK was partially supported by the Additional Funding Programme for Mathematical Sciences, delivered by EPSRC (EP/V521917/1) and the Heilbronn Institute for Mathematical Research. 

\section{The parametrization strategy}\label{sec:parametrization}

To introduce our strategy for producing algebraic points on $C$, we begin with a concrete example, to which we will return throughout.

\begin{Example}
    \label{ex:running_example_I_S6point}
    Consider the trigonal superelliptic curve
    \[C \colon y^3 = x^4 + 1,\]
    with $m=3$ and $d=4$. Suppose $\alpha$ is a root of $K = \Q[t]/F(t)$ where $F$ is the irreducible sextic 
    \[F = -t^6 + 3t^5 - 2t^4 + t^3 + 1\]
    which has $\Gal(\wt{K}/\Q) \simeq S_6$. A straightforward calculation reveals that 
    \[ \alpha^4 + 1 - (\alpha^2-\alpha)^3 = F(\alpha) = 0,\]
    Thus we have an $S_6$-sextic point $(\alpha, \alpha^2-\alpha) \in C(K)$. 
\end{Example}

In general, to produce algebraic points on $C$, our strategy is to parameterize the coordinates $x$ and $y$ as rational functions in an auxiliary variable $t$. Geometrically, this corresponds to producing rational curves and intersecting with $C$. Explicitly, we set
\[x(t) = \frac{\gamma(t)}{\eta(t)} \quad \text{and} \quad y(t) = \frac{g(t)}{h(t)}.\]
Substituting into the equation for $C$, given by \eqref{eq:sec}, and clearing denominators, we obtain the polynomial equation
\begin{equation}\label{eq:param}
F_{g,h,\gamma,\eta}(t) = h(t)^m \Big(c_d\gamma(t)^d + c_{d-1}\gamma(t)^{d-1}\eta(t) + \cdots + c_1\gamma(t)\eta(t)^{d-1} + c_0\eta(t)^d \Big) - g(t)^m\eta(t)^d = 0.\end{equation}
Suppose $g,h,\gamma, \eta$ are chosen in $\Z[x]$ such that $F_{g,h,\gamma, \eta}(t)$ is irreducible with some root $\alpha$. Then 
\[P = (x(\alpha), y(\alpha)) = \left(\frac{\gamma(\alpha)}{\eta(\alpha)}, \frac{g(\alpha)}{h(\alpha)}\right)\]
is a point on $C$ defined over the field $\Q(\alpha)$, and $\Q(\alpha)$ is the field generated by $P$. Given a degree $n$, our approach is to count how many ways we can choose $g,h,\gamma, \eta$ such that $F_{g,h,\gamma,\eta}$ is degree $n$, irreducible, and when possible, has Galois group $S_n$.

\begin{Example}\label{ex:running_example_II_S14point}
    Returning to our running example,
    \[C \colon y^3 = x^4 + 1,\]
    we see that the sextic point from Example \ref{ex:running_example_I_S6point} arose from this construction, taking $g(t) = t^2-t$, $h(t) = 1$, $\gamma(t) = t$, $\eta(t) = 1$.

    For another example of higher degree, suppose $n=14$: here we might take $\deg g = 4$, $\deg h = 2$, $\deg \gamma = 2$, and $\eta = 1$. Specifically setting $g = t^4 + t + 1$, $h= t^2$, $\gamma = t^2$, we obtain
    \[F = (t^2)^3f(t^2) - (t^4+t+1)^3=t^{14} - t^{12} - 3t^9 - 3t^8 - 2t^6 - 6t^5 - 3t^4 - t^3 - 3t^2 - 3t - 1\]
    which is irreducible with Galois group $S_{14}$. If $\alpha$ is a root of $F$, we have a degree 14 point $(\alpha^2, \frac{\alpha^4 + \alpha + 1}{\alpha^2}) \in C(K)$ for $K = \Q[t]/F(t)$.
\end{Example}

Generally, the degree of $F_{g,h,\gamma,\eta}$ is the maximum of $m(\deg h)+d(\deg \gamma)$ and $m(\deg g) + d(\deg \eta)$, both of which are multiples of $\gcd(m,d)$. Since we will eventually count the number of such parameterizations, we want to choose $g,h,\gamma,\eta$ so the sum of their degrees is as large as possible, giving us the most degrees of freedom to count. Recall that in this paper, we have assumed $m \leq d$, so this sum of degrees will be maximized by letting $\deg g$ and $\deg h$ be large, while keeping those of $\gamma$ and $\eta$ small. To that end, we simply take $\eta = 1$ and suppress the notation by writing $F_{g,h,\gamma}$ for the remainder of this paper. However, in the general case, namely if $m > d$, it would be useful to take $\eta$ to be nonconstant.

We observe that when $n$ is a sufficiently large multiple of $\gcd(m,d)$, we can always choose the degrees of $g$, $h$, and $\gamma$ to make the polynomial \eqref{eq:param} have degree $n$ in general. This is done by using $\deg \gamma$ to control the residue class of $n$ modulo $m$ if necessary, and letting $\deg g, \deg h$ be as large as possible. It remains to determine how large $n$ must be for such degrees to exist. It is clear that we must have at least $n \geq d$ by looking at the minimum degree of $F_{g,h,\gamma}$. To give a more precise answer we recall the classical definition of the Frobenius number, with a straightforward generalization to integers that are not coprime.

\begin{Definition}[Frobenius number]
	Given natural numbers $a,b$ with $\gcd(a,b) = 1$, the \textbf{Frobenius number}, denoted $\Frob(a,b)$, is the largest natural number which is not a linear combination $ax + by$ where $x,y \geq 0$. 

    When $\gcd(a,b) \neq 1$, we take $\Frob(a,b)$ to be the largest multiple of $\gcd(a,b)$ that is not a linear combination $ax + by$ for $x,y \geq 0$.
\end{Definition}

For coprime integers $a,b$, the Frobenius number is given by $\Frob(a,b) = ab - a - b$. Recognizing that for any natural numbers $a,b$ we have 
\[\Frob(a,b) / \gcd(a,b) = \Frob\left(\frac{a}{\gcd(a,b)}, \frac{b}{\gcd(a,b)}\right),\] 
we have $\Frob(a,b) = \lcm(a,b) - a - b$.

For any $n \geq \max(d, \Frob(m,d) +1)$ we can manipulate the degrees of $g$, $h$, and $\gamma$ such that $\deg F_{g,h,\gamma} = n$ in \eqref{eq:param}. Moreover, this is sharp in the sense that \eqref{eq:param} will not take degrees $n < d$ or $n = \Frob(m,d)$. We conclude this section by summarizing our discussion in the following  proposition.

\begin{Proposition}\label{prop:param}
Let $C$ be given by \eqref{eq:sec} with $m \leq d$. For all degrees $n \geq \max(d, \Frob(m,d) + 1)$ such that $\gcd(m,d) \mid n$, there exist $g,h,\gamma, \eta$ such that $F_{g,h,\gamma,\eta}(t)$ given in \eqref{eq:param} has degree $n$.

Explicitly, we can assume $\eta = 1$ and take $g,h,\gamma$ to have the degrees given below:
\begin{align}\label{eq:param_deg1}
\nonumber	\deg g &= n/m\\
	\deg h &= \lfloor (n-d)/m \rfloor & \text{when } m \mid n\\
\nonumber	\deg \gamma &= 1
\end{align}
and
\begin{align}\label{eq:param_deg2}
\nonumber	\deg g &= \lfloor n/m\rfloor\\
	\deg h &= (n-rd)/m & \text{when } m \nmid n\\
\nonumber	\deg \gamma &= r
\end{align}
where $r > 0$ is the minimal integer such that $n \equiv rd \pmod{m}$.
\end{Proposition}

Notice that the choices above accomplish our goals of maximizing the total degrees of freedom by letting $g,h$ have the largest possible degree, while $\deg \gamma$ is kept small, with $1 \leq r < m$.

\begin{Example}\label{ex:running_example_III_higherdeg}
    Let us return again to our running example curve
    \[C \colon y^3 = x^4 + 1.\]
    We have $\Frob(3,4)=6$, and $\gcd(m,d)=1$, so $n_0 = 6$. Thus for all $n \geq 6$, Proposition \ref{prop:param} suggests a family of candidate polynomials $F(t) \in \Q(\a, \b, \balpha)[t]$ for producing degree $n$ points on $C$.

    In the cases of $n=6,14$ addressed above in Examples \ref{ex:running_example_I_S6point} and \ref{ex:running_example_II_S14point}, we have seen that $F$ does indeed have irreducible specializations. Applying Hilbert's irreducibility theorem (see Lemma \ref{thm:hilbert_irred}), we find that $F$ almost always specializes to irreducible polynomials of degree $6$ with Galois group $S_6$ (respectively degree 14 with Galois group $S_{14}$). In this cases, we are ready to count these specializations and the fields they produce.
\end{Example}

To make this procedure work in general, we need to show that for any curve $C$ and sufficiently large degree $n$ divisible by $\gcd(m,d)$, the polynomial family $F_{g,h,\gamma}(t)$ of Proposition \ref{prop:param} is irreducible; this is proved in Proposition \ref{prop:param_irred_Sn}. For Theorem \ref{thm:m=3_bound}, we also need to show this family has Galois group $S_n$ when $m=3$ and $n \equiv 2,4 \pmod{6}$ is sufficiently large; this is accomplished in Proposition \ref{prop:m=3_irr_sn}. A key tool used in the proof of these results is the Newton polygon, which we introduce in the following section. The reader willing to grant that $F_{g,h,\gamma}(t)$ is irreducible may skip to Section \ref{sec:lower_bounds}, where we count the number of fields produced by specializations.

\section{Newton polygons}
\label{sec:galois_newton}

We now introduce the Newton polygon, which associates to a polynomial over $\Q_p$ a diagram of line segments containing data about valuations of roots. We will use this to show our polynomials are irreducible and identify cycles in the Galois group. Let $p$ be a prime, $\Q_p$ the field of $p$-adic numbers, and $F(t) \in \Q_p[t]$ a polynomial. 

\begin{Definition}[Newton polygon]
\label{def:newton_polygon}
	With the notation above, let $F(t)$ be given by $F(t) = \sum_{i=0}^n k_it^i$. The \textbf{$p$-adic Newton polygon of $F$} is the lower convex hull of the set
	\[\set{(i, v_p(k_i)) \in \R^2 \mid 0 \leq i \leq n},\]
	where $v_p$ denotes the $p$-adic valuation, and we set $v_p(0) = \infty$ by convention. We will denote the Newton polygon of $F$ by $\NP_{\Q_p}(F)$, or simply by $\NP(F)$ when it will not create confusion.
\end{Definition}

A good reference for the theory of Newton polygons is \cite[II.6]{Neukirch}. In particular, the Newton polygon $\NP(F)$ can be split up into segments of distinct slopes $s_j$, and if the $j$-th segment has length $\ell_j$, then $F(t)$ has $\ell_j$ roots of valuation $-s_j$. This key fact leads to the following lemmas, proven in \cite{keyes}.

\begin{Lemma}[{see \cite[Lemma 2.6]{keyes}}]
\label{lem:NP_factorization}
	Suppose $\NP_{\Q_p}(F)$ has a segment of length $\ell$ and slope $s$, and no other segments of this slope (i.e.\ consider the \textit{entire} segment of slope $s$). Then $F$ factors as $F = F_0F_1$ over $\Q_p$, where $\deg F_0 = \ell$ and the roots of $F_0$ have $p$-adic valuation $-s$.
	
	Moreover, if $s = r/\ell$ has reduced fraction form $r'/\ell'$ then all irreducible factors of $F_0$ over $\Q_p$ have degree divisible by $\ell'$. In particular, if $\gcd(r,\ell) = 1$ then the $F_0$ produced above is irreducible.
\end{Lemma}

\begin{Lemma}[{see \cite[Lemma 2.7]{keyes}}]
\label{lem:NP_cycle}
	Suppose $F(t) \in \Q[t]$, $p > \deg F$, and $\NP_{\Q_p}(F)$ has a segment of length $\ell$ and slope $r/\ell$ with $\gcd(r,\ell) = 1$. Let $F = F_0F_1$ be the factorization of Lemma \ref{lem:NP_factorization}, for which $F_0$ is irreducible. If $\ell$ is pairwise coprime to the degrees of the irreducible factors of $F_1$ over $\Q_p$, then $\Gal(F/\Q)$ contains an $\ell$-cycle.
\end{Lemma}

When $F \in \Z[t]$, Lemma \ref{lem:NP_factorization} may be used to deduce that $F$ is irreducible by applying it at one or more primes. In this case, the Galois group of $F$ is a transitive subgroup of $S_n$. Lemma \ref{lem:NP_cycle} may be used to show certain cycles are present in $G$, which sometimes suffices to determine $G = S_n$. We will make use of the following standard result.

\begin{Lemma}
\label{lem:genset}
	Let $G \subseteq S_n$ be a permutation subgroup acting on the set $\lbrace 1, \ldots, n\rbrace$. Suppose the action of $G$ is transitive and that $G$ contains a transposition. Then if $G$ contains a $p$-cycle for a prime $p > n/2$, we have $G = S_n$.
\end{Lemma}

\begin{proof}
	This is a standard group theory exercise. For a short proof see \cite[Proposition 2.4]{keyes}.
\end{proof}

See also \cite[\S 2.1.3]{keyes_thesis} for proofs of these lemmas, along with further discussion and examples of how they can be used to elucidate the Galois group of polynomials $F(t)$. 

\section{Irreducibility of $F(t)$}\label{sec:irred_sn}

Let $C$ be a superelliptic curve with exponent $m$ and defining polynomial $f(x)$, as in \eqref{eq:sec}. As in Proposition \ref{prop:param}, given any $n \geq \max(d, \Frob(m,d) + 1)$ such that $\gcd(m,d) \mid n$, there exist choices of degrees \eqref{eq:param_deg1} or \eqref{eq:param_deg2} for $g,h, \gamma$ such that the polynomial $F_{g,h,\gamma}(t)$ given in \eqref{eq:param} has degree $n$ in general. Writing
\begin{align*}
	g(t) = \sum_{i=1}^{\deg g} a_i t^i, \quad\quad	h(t) = \sum_{j=1}^{\deg h} b_j t^j,  \quad\quad	\gamma(t) = \sum_{\ell=1}^{\deg \gamma} \alpha_\ell t^\ell,
\end{align*}
we can view $F_{g,h,\gamma}(t)$ as a degree $n$ polynomial $F(\a, \b, \boldsymbol{\alpha},t) \in \Q(\a, \b, \boldsymbol{\alpha})[t]$. Here $\a$ indicates the tuple of indeterminates $(a_0, \ldots, a_{\deg g})$, and similarly for $\b$ and $\balpha$. For simplicity, since we have fixed the curve $C$ and degree $n$, we will denote this polynomial family by $F \in \Q(\a,\b, \boldsymbol{\alpha})[t]$, and denote a rational specialization by $F_{\a_0, \b_0, \boldsymbol{\alpha}_0} \in \Q[t]$, where $\a_0 \in \Q^{\deg g + 1}, \b_0 \in \Q^{\deg h + 1}, \boldsymbol{\alpha}_0 \in \Q^{\deg \gamma + 1}$. 

\newcommand{\y}{\boldsymbol{y}}

Since $F$ is degree $n$, almost all specializations $F_{\a_0, \b_0, \boldsymbol{\alpha}_0}$ have degree $n$. We recall Hilbert's irreducibility theorem, which states that the irreducibility of the polynomial family carries over to almost all specializations. We state this for a general polynomial $F(\y,t) \in \Q(\y)[t]$ where $\y$ is some tuple of indeterminates.

\begin{Lemma}[Hilbert's irreducibility theorem]
\label{thm:hilbert_irred}
	Let $F(\y, t) \in \Q(\y)[t]$ be irreducible with Galois group $G$. Then for 100\% of specializations $\y_0$, we have $F(\y_0, t) \in \Q[t]$ is irreducible with Galois group $G_0 \simeq G$.
\end{Lemma}

Using the Newton polygons from the previous section, our aim is to study the factorizations over $\Q_p$ of integral specializations $F_{\a_0, \b_0, \boldsymbol{\alpha}_0}$ to show that $F$ is irreducible over $\Q(\a, \b, \boldsymbol{\alpha})$. Lemma \ref{thm:hilbert_irred} then implies that almost all specializations $F_{\a_0, \b_0, \boldsymbol{\alpha}_0}$ are also irreducible. The remainder of this section is devoted to proving the following proposition, which makes this precise when $n$ is sufficiently large.

\begin{Proposition}
\label{prop:param_irred_Sn}
	Fix an integer $m \geq 2$, a degree $d\geq m$ polynomial $f(x) \in \Z[x]$ satisfying \eqref{eq:f_not_power}, and an integer $n \geq n_0 = \max(d, \Frob(m,d)+1)$ such that $\gcd(m,d) \mid n$. The degree $n$ polynomial family $F(t) \in \Q(\a,\b,\boldsymbol{\alpha})[t]$, given in \eqref{eq:param} with degrees \eqref{eq:param_deg1} if $m \mid n$ or \eqref{eq:param_deg2} if $m \nmid n$, is irreducible whenever one of the following is satisfied:
    \begin{enumerate}[label = (\roman*)]
        \item $m \mid n$;
        \item $\gcd(n,m) = 1$;
        \item $n \geq \frac{m^2(m-1)}{2}$.
    \end{enumerate}	
	Moreover, when $F$ is irreducible, 100\% of specializations $F_{\a_0, \b_0, \balpha_0}$ are irreducible of degree $n$ over $\Q$.
\end{Proposition}

The second statement follows from the first by Hilbert irreducibility, Lemma \ref{thm:hilbert_irred}. Before proving the first statement, we need some elementary results. 

\begin{Lemma}\label{lem:primes_dividing_f(x)}
	Let $f(x) \in \Z[x]$ be a nonconstant polynomial. Then there exist infinitely many primes $p$ such that $p \mid f(x_0)$ for some integer $x_0$. 
 
    Moreover, if $f(x)$ has an irreducible factor $f_0$ appearing with multiplicity $e \geq 1$, then there exist infinitely many primes $p$ such that $v_p(f(x_0)) = e$ for some integer $x_0$.
\end{Lemma}

\begin{proof}
	A prime $p$ divides $f(x_0)$ for some integer $x_0$ if and only if the reduction of $f$ modulo $p$ has a root. The set of such primes certainly contains the primes for which the reduction of $f$ modulo $p$ splits completely. These are precisely the primes which split completely in the splitting field of $f$. By the Chebotarev Density Theorem, this is an infinite set, proving the first claim.

    We prove the second claim first for squarefree $f$, in which case $e=1$. Suppose $p \nmid \Disc f$ and that $p \mid f(x_0)$ for some integer $x_0$, the existence of which is guaranteed by the first claim. Consider $f(x_0 + p)$,
	\[f(x_0 + p) \equiv f(x_0) + f'(x_0)p \pmod{p^2}.\]
	If $p^2 \mid f(x_0)$ and $p^2 \mid f(x_0 + p)$ then we must have $p \mid f'(x_0)$. However, this implies $x_0$ is a double root of $f(x)$ mod $p$, contradicting $p \nmid \Disc f$. Thus we conclude one of $f(x_0)$ and $f(x_0 + p)$ is divisible by $p$ exactly once, hence the set of primes dividing $f(x_0)$ exactly once for some $x_0$ is infinite.

    Let us now remove the assumption that $f$ is squarefree and write its irreducible factorization $f = \prod_{i \geq 0} f_i^{e_i}$ for $e_i \geq 1$. Set $g = \prod_{i \geq 0} f_i$, which is squarefree by construction. By the first claim of the lemma, there exists a prime $p > \Disc(g)$ and an integer $x_0$ such that $v_p(f_0(x_0)) = 1$. Since $p > \Disc(g)$, $p \nmid \prod_{i > 0} f_i(x_0)$. Upon returning to $f$, we have $v_p(f_0(x_0)^{e_0}) = e_0$ and $p \nmid f_i(x_0)^{e_i}$, completing the proof of the second claim. 
\end{proof}

\begin{Lemma}\label{lem:divisibility_fraction}
    Let $n, u,v,w$ be positive integers such that $u,v \mid n$. If $\frac{n}{u}, \frac{n}{v} \mid w$ then $\frac{n}{\gcd(u,v)} \mid w$.
\end{Lemma}

\begin{proof}
    Omitted.
\end{proof}

We are now ready to prove Proposition \ref{prop:param_irred_Sn}.

\begin{proof}[Proof of Proposition \ref{prop:param_irred_Sn}]
    For the first statement, we look separately at the cases of $m \mid n$, $\gcd(m,n) = 1$, and $1 < \gcd(m,n) < m$.

    \subsubsection*{Case (i): $m \mid n$} Our goal is to exhibit specializations with incompatible $p$-adic factorizations for several primes $p$, arguing via Newton polygons and Lemma \ref{lem:NP_factorization}.

    Fix a prime $p$ such that $p \nmid c_i$ for all $i$. Consider an integral specialization $\a_0, \b_0, \boldsymbol{\alpha}_0$ satisfying
	\begin{align}
	\label{eq:irr1}
		v_p(a_0) &= 1 \\
\nonumber		v_p(a_i) &\geq 1 \text{ for } 0 < i < n/m\\
\nonumber		v_p(a_{n/m}) &= 0 \\
\nonumber		v_p(b_j) &\geq 1 \text{ for } 0 \leq j \leq (n-d)/m,
	\end{align}
	with no restrictions on $\alpha_0, \alpha_1$. 	We end up with the Newton polygon featured below in Figure \ref{fig:sec_irr1}.	
\begin{figure}[ht]
\centering
\caption{$\NP_{\mathbb{Q}_{p}}(F_{\a_0,\b_0,\boldsymbol{\alpha_0}})$ with one segment of slope $-m/n$}
\label{fig:sec_irr1}
\vspace{1ex}

\begin{tikzpicture}[scale=0.75]
	\draw[->, thick] (-0.2, 0) -- (10.5, 0);
	\draw[->, thick] (0, -0.2) -- (0, 2.5);
	
	\filldraw[black] (0, 2) circle (2pt) node[left] {$(0, m)$};
	\filldraw[black] (10,0) circle (2pt) node[below] {$(n, 0)$};
	
	\draw[black] (0,2) -- (10,0);
\end{tikzpicture}
\end{figure}	
	
	In particular, since we have assumed $m \mid n$, we have that $m = \gcd(m,n)$. By Lemma \ref{lem:NP_factorization}, all irreducible factors of $F$ over $\Q_p$ must have degree divisible by $\frac{n}{m}$.
	
	Consider now an alternative specialization. By Lemma \ref{lem:primes_dividing_f(x)} there are infinitely many primes $p$ such that for some $\alpha_0 \in \Z$ we have $f(\alpha_0)$ is divisible by $p$ exactly $e$ times, where $e=e_i$ is the multiplicity of an irreducible factor of $f(x) = \prod_i f_i(x)^{e_i}$. Note that we may not be able to enforce $e=1$ since $f$ need not be squarefree, nor even have an irreducible factor of multiplicity one.
	
    Choose some such $p$ and $\alpha_0$ such that $p\nmid c_i$ for all $i$. For any $\gamma(t) \equiv p^et + \alpha_0 \pmod{p^{e+1}}$, we have that $p^e$ divides all coefficients of $f(\gamma(t))$, and in particular $p^e \ || \ f(\alpha_0)$. Consider now a specialization satisfying
 \begin{align}
	\label{eq:irr2}
		v_p(a_i) &\geq e \text{ for } 0 \leq i < n/m \\
\nonumber		v_p(a_{n/m}) &= 0 \\
\nonumber		v_p(b_0) &= 0 \\
\nonumber		v_p(b_j) &\geq 0 \text{ for } 0 \leq j \leq (n-d)/m.
	\end{align}
	This ensures that $p^e$ exactly divides the constant term of $F$ and all other terms except the leading term, yielding the Newton polygon below in Figure \ref{fig:sec_irr2}.
	\begin{figure}[ht]
\centering
\caption{$\NP_{\mathbb{Q}_{p}}(F_{\a_0,\b_0,\boldsymbol{\alpha_0}})$ with one segment of slope $-e/n$}
\label{fig:sec_irr2}
\vspace{1ex}

\begin{tikzpicture}[scale=0.75]
	\draw[->, thick] (-0.2, 0) -- (10.5, 0);
	\draw[->, thick] (0, -0.2) -- (0, 1.5);
	
	\filldraw[black] (0, 1) circle (2pt) node[left] {$(0, e)$};
	\filldraw[black] (10,0) circle (2pt) node[below] {$(n, 0)$};
	
	\draw[black] (0,1) -- (10,0);
\end{tikzpicture}
\end{figure}
	
	As earlier, Lemma \ref{lem:NP_factorization} implies that the an irreducible factor of $F$ over $\Q_p$ must have degree a multiple of $\frac{n}{\gcd(n,e)}$. 
	
	Let $F_0$ be an irreducible factor of $F$. We have seen that 
	\[\frac{n}{m} \ \Big| \ \deg F_0 \quad \text{and} \quad \frac{n}{\gcd(n,e_i)} \ \Big| \ \deg F_0 \text{ for all } e_i.\]
 Applying Lemma \ref{lem:divisibility_fraction} to each $e_i$ (with $u=m,\ v= \gcd(n, e_1, \ldots, e_{i-1}),\ w = \deg F_0$), we find
	\[\frac{n}{\gcd(n,m,e_i)} = \frac{n}{1} \ \Big| \ \deg F_0\]
	by our assumption \eqref{eq:f_not_power}. Hence $F$ is irreducible over $\Q$.
    
    \subsubsection*{Case (ii): $\gcd(m,n) = 1$} In this case, the degrees of $g(t), h(t)$, and $\gamma(t)$ are given in \eqref{eq:param_deg2}. Let $p$ be a prime such that $p \nmid c_0, c_d$. We consider a specialization $\a_0, \b_0, \boldsymbol{\alpha}_0$ satisfying the following requirements.
	\begin{align}\label{eq:case_a_transitive}
		v_p(a_0) &= 0 \\
		\nonumber v_p(a_i) & \geq m \text{ for all } i > 0\\
		\nonumber v_p(b_j) & \geq m \text{ for all } j < (n-rd)/m\\
		\nonumber v_p(b_{(n-rd)/m}) &= 1\\
		\nonumber v_p(\alpha_{\ell}) &\geq 0 \text{ for all } \ell < r\\
		\nonumber v_p(\alpha_r) &= 0.
	\end{align}
The choices in \eqref{eq:case_a_transitive} ensure that all but the constant term of $F_{\a_0, \b_0, \balpha_0}$ have $p$-adic valuation at least $m$, while the leading term of is $b_{(n-rd)/m}^mc_d\alpha_r^d$, which has $p$-adic valuation exactly $m$. The constant term, $b_0^mf(\alpha_0) - a_0^m$ has valuation 0. This produces the $p$-adic Newton polygon below in Figure \ref{fig:case_a_transitive}.

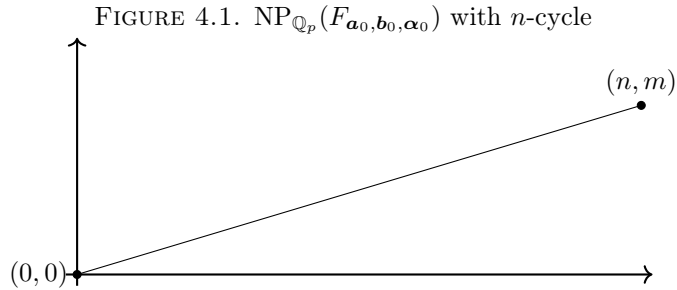
\begin{figure}[H]
\centering
\caption{$\NP_{\Q_p}(F_{\a_0,\b_0,\balpha_0})$ with $n$-cycle}
\label{fig:case_a_transitive}

\begin{tikzpicture}[scale=0.75]
	\draw[->, thick] (-0.2, 0) -- (10.5, 0);
	\draw[->, thick] (0, -0.2) -- (0, 2.5);
	
	\filldraw[black] (0, 0) circle (2pt) node[left] {$(0, 0)$};
	\filldraw[black] (10, 2) circle (2pt) node[right] {$(n, m)$};
	
	\draw[black] (0,0) -- (10,2);
	
\end{tikzpicture}

\end{figure}
\noindent Since the polygon has exactly one segment of slope $m/n$ with $\gcd(m,n) = 1$, Lemma \ref{lem:NP_factorization} implies that $F_{\a_0, \b_0, \balpha_0}$ is irreducible. Hence $F$ must be irreducible.

    \subsubsection*{Case (iii): $1 < \gcd(m,n) < m$} 

    For this final case, we blend the strategies of the previous cases: if $F_0$ is an irreducible factor of $F$ of degree $d_0$, then the approach from Case (ii) shows that $\frac{n}{\gcd(n,m)} \mid d_0$, while that of Case (i) shows that $d_0$ is close to a multiple of $\frac{m\lfloor \frac{n}{m}\rfloor}{\gcd(m\lfloor \frac{n}{m} \rfloor, e)}$, where $e$ is the multiplicity of an irreducible factor of $f$. The additional hypothesis that $n \geq \frac{m^2(m-1)}{2}$ allows us to apply Lemma \ref{lem:irreducible_conspiracy}, a technical intermediate that allows us to conclude $d_0$ is a multiple of $\frac{n}{\gcd(n,m,e)}$.    
    
    More precisely, the degrees of $g(t), h(t)$, and $\gamma(t)$ are given by \eqref{eq:param_deg2}. Choosing some prime $p$ not dividing $c_0$ or $c_d$, the restrictions \eqref{eq:case_a_transitive} produce the $p$-adic Newton polygon in Figure \ref{fig:case_a_transitive}. However, since $\gcd(m,n) >1$, we cannot conclude right away that $F_{\a_0, \b_0, \balpha_0}$ is irreducible. Since the segment has slope $m/n = \frac{m/\gcd(m,n)}{n/\gcd(m,n)}$, Lemma \ref{lem:NP_factorization} gives that any irreducible factors of $F_{\a_0, \b_0, \balpha_0}$ has degree divisible by $n/\gcd(m,n)$. In particular, this means that any irreducible components of $F$ must also have degree divisible by $n/\gcd(m,n)$.

    Consider now an irreducible factor $f_0$ of $f$ appearing with multiplicity $e$ in the irreducible factorization. By Lemma \ref{lem:primes_dividing_f(x)}, there exists another prime $p$ and an integer $\alpha_0$ such that $v_p(f(\alpha_0)) = e$. Take the restrictions \eqref{eq:irr2}, noting that we must set $v_p(a_{\lfloor n/m\rfloor}) = 0$, since $m \nmid n$. This produces a Newton polygon similar to that of Figure \ref{fig:sec_irr2}, but with a segment of length $m\lfloor \frac{n}{m} \rfloor$ and slope $\frac{-e}{m\lfloor \frac{n}{m} \rfloor}$. This produces the polygon shown below in Figure \ref{fig:case_c_long_segment}.

    \begin{figure}[ht]
    \caption{$\NP_{\mathbb{Q}_{p}}(F_{\a_0,\b_0,\boldsymbol{\alpha_0}})$ with one segment of slope $\frac{-e}{m\lfloor n/m \rfloor}$}
    \label{fig:case_c_long_segment}
    \vspace{1ex}

    \begin{tikzpicture}[scale=0.75]
	   \draw[->, thick] (-0.2, 0) -- (10.5, 0);
	   \draw[->, thick] (0, -0.2) -- (0, 2.5);
	
	   \filldraw[black] (0, 1) circle (2pt) node[left] {$(0, e)$};
	   \filldraw[black] (8,0) circle (2pt) node[below] {$(m\lfloor\frac{n}{m} \rfloor, 0)$};
          \filldraw[black] (10,2) circle (2pt) node[right] {$(n,\geq re)$};
	
	   \draw[black] (0,1) -- (8,0);
	   \draw[black, dashed] (10,2) -- (8,0);
    \end{tikzpicture}
    \end{figure}

    Note that the rightmost segment is not specified by our restrictions on the coefficients, indicated by the dashed line. We will not need to assume anything about this segment for this proof.

    From the Newton polygon in Figure \ref{fig:case_c_long_segment}, we deduce that in order to be compatible with the factorization of $F$ over $\Q_p$, an irreducible factor of $F$ must have degree
    \begin{equation}\label{eq:almost_multiple}
        k \frac{m\lfloor\frac{n}{m} \rfloor}{\gcd(m\lfloor\frac{n}{m} \rfloor, e)} + \ell
    \end{equation}
    for integers $0 \leq k \leq \gcd(m\lfloor\frac{n}{m} \rfloor,e)$ and $0 \leq \ell \leq n - m\lfloor\frac{n}{m} \rfloor$. When $n$ is sufficiently large, multiples of $\frac{n}{\gcd(n,m)}$ of the form \eqref{eq:almost_multiple} are in fact multiples of $\frac{n}{\gcd(n,m,e)}$; in Lemma \ref{lem:irreducible_conspiracy} below, we show that $n \geq \frac{m^2(m-1)}{2}$ is sufficient for this purpose.

    We may now apply Lemma \ref{lem:divisibility_fraction} as in Case (i) for all $e = e_i$ appearing in the irreducible factorization of $f$. This shows that the degree of any irreducible factor of $F$ is divisible by $\frac{n}{\gcd(n,m,e_i)} = n$, hence $F$ itself is irreducible.

    \bigskip \noindent This concludes the proof of the first statement. The second follows directly by Hilbert's irreducibility theorem, Lemma \ref{thm:hilbert_irred}, completing the proof of the proposition.
\end{proof}

\begin{Lemma}\label{lem:irreducible_conspiracy}
    Let $m, n, e$ be positive integers satisfying $1 < \gcd(m,n) < m$, $e < m$ and 
    \begin{equation}
	\label{eq:reqs_on_n}
            n \geq \frac{m^2(m-1)}{2}.
    \end{equation}
    Suppose that $d_0 \leq n$ is a positive integer multiple of $\frac{n}{\gcd(n,m)}$ satisfying
    \[d_0 = k\frac{m\lfloor \frac{n}{m} \rfloor}{\gcd(e,m\lfloor \frac{n}{m} \rfloor)} + \ell\]
    for nonnegative integers $k$, $\ell \leq n - m\lfloor \frac{n}{m} \rfloor$. Then we have
    \[\frac{n}{\gcd(n,m,e)} \mid d_0.\]
\end{Lemma}

To prove Lemma \ref{lem:irreducible_conspiracy}, we will use the following fact.

\begin{Lemma}\label{lem:divisor_interval}
    Fix an integer $n \geq 1$ and a real number $w > 0$. There is at most one divisor of $n$ in the open interval
    \[\left( \frac{n}{\lfloor w \rfloor + 1}, \frac{n}{\lceil w \rceil - 1}\right).\]
    Moreover, such a divisor exists if and only if $w$ is an integer divisor of $n$.
\end{Lemma}

\begin{proof}
    Suppose $w'$ is an integer dividing $n$ and $\frac{n}{w'}$ is contained in the open interval. Thus we have
    \[\lceil w \rceil - 1 < w' < \lfloor w \rfloor + 1.\]
    If $w$ is not an integer, we have $w' \in (\lfloor w \rfloor, \lfloor w \rfloor + 1)$, which contradicts $w'$ being an integer. If $w$ is an integer then $w' \in (w-1, w+1)$, forcing $w' = w$, in which case $\frac{n}{w}$ is the only possible divisor in the given interval.
\end{proof}

\begin{proof}[Proof of Lemma \ref{lem:irreducible_conspiracy}]
    If $m=1$ the statement is vacuously true, so assume $m \geq 2$. If $k=0$ then $d_0 = \ell \leq n-m\lfloor \frac{n}{m} \rfloor < m$, contradicting \eqref{eq:reqs_on_n}; moving forward we assume $k \geq 1$. 
    We write $d_0 = \frac{un}{\gcd(n,m)}$ for a positive integer $u$. To ease notation, we set $e' = \gcd(e,m\lfloor \frac{n}{m} \rfloor)$.
    
    We look to apply Lemma \ref{lem:divisor_interval} as follows: if we can show that $d_0$ is contained in the open interval
	\[\left( \frac{un}{\lfloor \frac{ue'}{k} \rfloor +1}, \frac{un}{\lceil \frac{ue'}{k} \rceil -1} \right)\]
	then we must have $\frac{ue'}{k}$ is integral and also a divisor of $n$. Moreover, we find that $d_0 = \frac{nk}{e'}$, which is divisible by $\frac{n}{\gcd(n, e')}$. Lemma \ref{lem:divisibility_fraction} then gives $\frac{n}{\gcd(n,m,e)} \mid d_0$, as desired.

    To show that the hypothesis \eqref{eq:reqs_on_n} ensures $d_0 < \frac{un}{\lceil \frac{ue'}{k} \rceil - 1}$, it suffices to show
    \[un - \left( \left\lceil \frac{ue'}{k} \right\rceil - 1 \right) d_0 = \frac{nk}{e'} \left(\frac{ue'}{k} - \left\lceil \frac{ue'}{k} \right\rceil + 1 \right) + \left(\left\lceil \frac{ue'}{k} \right\rceil - 1\right) \left( \ell - \frac{k(n-m\lfloor\frac{n}{m} \rfloor)}{e'}\right) > 0.\]
    Since $\frac{ue'}{k} - \lceil \frac{ue'}{k} \rceil + 1 \geq \frac{1}{k}$, we have that the first term is bounded below by $\frac{n}{e'} \geq \frac{n}{m}$. For the second term, we recall $\ell \geq 0$ and $\lceil \frac{ue'}{k} \rceil - 1 < \frac{ue'}{k}$
    \begin{align*}
        \left(\left\lceil \frac{ue'}{k} \right\rceil - 1\right) \left( \ell - \frac{k(n-m\lfloor\frac{n}{m} \rfloor)}{e'}\right) & > - u\left(n - m \left\lfloor \frac{n}{m} \right\rfloor\right) \\
        & \geq -\frac{m(m-1)}{2}.
    \end{align*}
    For the final inequality, we have used $n - m\lfloor\frac{n}{m}\rfloor < m$ and $u = \frac{d_0\gcd(n,m)}{n} \leq \gcd(n,m) \leq \frac{m}{2}$, as a consequence of our hypotheses on $m$ and $n$.

    Combining these, we have that 
    \[un - \left( \left\lceil \frac{ue'}{k} \right\rceil - 1 \right) d_0 > \frac{n}{m} - \frac{m(m-1)}{2} \geq 0,\]
    where the second inequality is equivalent to the hypothesis \eqref{eq:reqs_on_n}. A similar argument shows that \eqref{eq:reqs_on_n} suffices to ensure $d_0 > \frac{un}{\lfloor \frac{ue'}{k} \rfloor + 1}$, completing the proof.	
\end{proof}

\begin{Example}\label{ex:running_example_IV_irr}
    Let us return to our running concrete example (see Examples \ref{ex:running_example_I_S6point}, \ref{ex:running_example_II_S14point}, \ref{ex:running_example_III_higherdeg}),
    \[C \colon y^3 = x^4 + 1.\]
    Since in this case $m=3$ is prime, for all $n \geq n_0 = 6$, either $3 \mid n$ or $\gcd(n,3)=1$ so the hypotheses (i) or (ii) of Proposition \ref{prop:param_irred_Sn} are always satisfied. Thus, the polynomial family $F(t)$ described in Proposition \ref{prop:param} is irreducible over $\Q(\a, \b, \balpha)[t]$, and thus 100\% of specializations $F_{\a_0, \b_0, \balpha_0}(t)$ are irreducible over $\Q$. This generalizes our observations in the previous examples from the special cases of $n=6,14$ to all $n \geq 6$.
\end{Example}

\begin{Example}
	To illustrate why a hypothesis like \eqref{eq:reqs_on_n} is needed in Case (iii), consider the case of $m=10$ and $n = 24$. We see that 24 is far too small to satisfy \eqref{eq:reqs_on_n}. If $F_0$ is an irreducible factor of $F$, by the Newton polygon construction shown in Figure \ref{fig:sec_irr1} and Lemma \ref{lem:NP_factorization}, $\deg F_0$ is divisible by $\frac{n}{\gcd(n,m)} = 12$, i.e.\ at worst, $F$ has two factors of equal degree.

    Suppose $f = f_1^5f_2^2$ for irreducible factors $f_1, f_2$. Note that if, for example, $\deg f_1 = 2$ and $\deg f_2 = 1$, we have $\deg f = 12$, so $n=24$ is indeed attainable through our parametrization strategy.
	
	Using the argument described in the proof above, we cannot rule out that $F$ is reducible. Arguing with $e_1 = 5$ we have that $\deg F_0 = 4k + \ell$ for integers $k, 0 \leq \ell \leq 4$. Of course, $12 = 4(3) + 0 = 4(2) + 4$, so this is consistent with $\deg F_0 = 12$. Looking instead at $e_2 = 2$, we have that $\deg F_0 = 10k + \ell$, which is also consistent with $\deg F_0 = 12$.

    When $n$ is sufficiently large, as made precise by Lemma \ref{lem:irreducible_conspiracy}, such coincidences cannot occur. The specializations we produce in the proof above have incompatible factorizations, allowing us to conclude $F$ is irreducible.
\end{Example}

\section{Galois groups in the $m=3$ case}\label{sec:m=3}

Since $F$ is irreducible, by taking specializations we obtain infinitely many degree $n$ number fields generated by points on $C$. It is natural to ask about the structure of their Galois groups, or more precisely to identify the group $G = \Gal(F/\Q(\a,\b,\balpha))$. Generically, we might expect the Galois groups of the polynomials produced by our parametrization to be as large as possible, i.e. $G=S_n$. This is known in the $m=2$ case \cite{LOT, keyes}.

We investigate the case of trigonal ($m=3$) superelliptic curves and give a positive result for certain degrees $n$, which is a key step in the proof of Theorem \ref{thm:m=3_bound}. 

\begin{Proposition}\label{prop:m=3_irr_sn}
    Fix $m=3$ and an integral polynomial $f(x)$ of degree $d > 3$ with $3 \nmid d$ satisfying \eqref{eq:f_not_power}. For all $n \geq \max(2d-2, 14)$ and $n \equiv 2,4 \pmod{6}$, the polynomial family $F(t)$ given in \eqref{eq:param_deg2} has Galois group $G = S_n$.
\end{Proposition}

By Proposition \ref{prop:param_irred_Sn}, we have that $F$ is irreducible, so $G \subseteq S_n$ is a transitive permutation subgroup. To show equality, we find certain cycle types in $G$ and apply Lemma \ref{lem:genset}. Along the way, we use the following corollary of Hilbert's irreducibility theorem.

\begin{Corollary}[see {\cite[Theorem 4.2]{LOT}}]
\label{cor:hilbert_irred}    
    Suppose $F(\y,t) \in \Q(\y)[t]$ is irreducible. If a permutation representation of $G_0 = \Gal(F(\y_0, t)/\Q)$ contains a given cycle type for a positive proportion of integral specializations $\y_0$, then $G$ contains an element of the same cycle type.
\end{Corollary}

Let us first show that $G$ contains a transposition. To do this, we use Lemma \ref{lem:NP_cycle}; the challenge is to find a prime $p$ and a specialization $\a_0, \b_0, \balpha_0$ depending only on the residue class modulo a power of $p$ such that $\NP_{\Q_p}(F_{\a_0, \b_0, \balpha_0})$ has one segment of length 2 and no irreducible factors of odd degree.

\begin{Lemma}\label{lem:m=3_transp}
    Fix $m=3$ and an integral polynomial $f(x)$ of degree $d > 3$ with $3 \nmid d$ satisfying \eqref{eq:f_not_power}. If $n \geq 2d-2$ and $n \equiv 2, 4 \pmod{6}$, then $G$ contains a transposition.
\end{Lemma}

\begin{proof}
    The condition that $n \geq 2d - 2$ is equivalent to the hypothesis $n \geq \max(d, \Frob(3,d)+1)$ from Proposition \ref{prop:param_irred_Sn}. Hypothesis (ii) of that proposition is also satisfied by our congruence conditions on $n$, so we have that $F$ is irreducible. 

    Let us look more closely at the low degree and leading terms of $F_{\a_0, \b_0, \balpha_0}(t)$:
\begin{align}\label{eq:low_deg_terms_F}
    F_{\a_0, \b_0, \balpha_0}(t) &= a_0^3 - b_0^3f(\alpha_0) \\
    \nonumber &+ \left(3a_0^2a_1 - b_0^3f'(\alpha_0)\alpha_1 - 3b_0^2b_1f(\alpha_0)\right)t \\
    \nonumber &+ \left(3a_0^2a_2 + 3a_0a_1^2 - b_0^3\left(\frac{f''(\alpha_0)}{2}\alpha_1^2 + f'(\alpha_0)\alpha_2 \right) - 3b_0^2b_1 f'(\alpha_0)\alpha_1 - (3b_0^2b_2 + 3b_0b_1^2)f(\alpha_0) \right)t^2\\
    \nonumber & \vdots \\
    \nonumber & + b_{(n-rd)/3}^3 c_d \alpha_r^d t^n.
\end{align}
(Note if $\ell > r= \deg \gamma$ then we take $\alpha_\ell = 0$ above.)

Suppose $f$ has at least one irreducible factor of multiplicity $e_i = 1$; equivalently, $f$ is not a perfect square. We will deal with the other case later. It follows from the same argument as the proof of Lemma \ref{lem:primes_dividing_f(x)} that there exists a prime $p\nmid c_d$ and an integer $\alpha_0$ such that $p^2 \ || \ f(\alpha_0)$ and $p \nmid f'(\alpha_0)$. Consider an integral specialization $\a_0, \b_0, \balpha_0$ satisfying the following conditions:
\begin{align}\label{eq:m=3_e=1_transp}
    v_p(a_0) &\geq 3 \\
    \nonumber v_p(a_1) &= 0 \\
    \nonumber v_p(a_i) & \geq 3 \text{ for all } 1 < i \leq \lfloor \frac{n}{3} \rfloor \\
    \nonumber v_p(b_j) &= 1 \text{ for all } 0 \leq j \leq (n-rd)/3 \\
    \nonumber v_p(\alpha_\ell) &= 0 \text{ for all } 0 < \ell \leq r.
\end{align}
These conditions ensure that each coefficient of $F$ is divisible by $p^3$, except for the $t^3$ term, which is not divisible by $p$. Inspecting the terms of $F$ in \eqref{eq:low_deg_terms_F} more carefully, we see that the constant, linear, and leading terms have valuations 5, 3, and 3, respectively, producing the Newton polygon below in Figure \ref{fig:n_odd_transp}
\begin{figure}[ht]
    \centering
    \caption{$\NP_{\mathbb{Q}_{p}}(F_{\a_0,\b_0,\boldsymbol{\alpha_0}})$ with transposition}
    \label{fig:n_odd_transp}
    \vspace{1ex}

    \begin{tikzpicture}[scale=0.75]
	   \draw[->, thick] (-0.2, 0) -- (10.5, 0);
	   \draw[->, thick] (0, -0.2) -- (0, 3);
	
	   \filldraw[black] (0, 2.5) circle (2pt) node[left] {$(0, 5)$};
	   \filldraw[black] (1, 1.5) circle (2pt) node[right] {$(1,3)$};
          \filldraw[black] (3, 0) circle (2pt) node[below] {$(3,0)$};
          \filldraw[black] (10, 1.5) circle (2pt) node[above] {$(n, 3)$};
	
	   \draw (0, 2.5) -- (1, 1.5);
	   \draw (3,0) -- (1, 1.5);
	   \draw (3,0) -- (10,1.5);
    \end{tikzpicture}
    \end{figure}
The rightmost segment has length $n-3$ and slope $\frac{3}{n-3}$. Since $n \equiv 2$ or $4 \pmod{6}$, we have both that $n-3$ is odd and $\gcd(3, n-3) = 1$. The polygon therefore satisfies the hypotheses of Lemma \ref{lem:NP_cycle}, producing a transposition (and in fact an $(n-3)$-cycle) in the Galois group $\Gal(F_{\a_0, \b_0, \balpha_0}/\Q)$. Since the restrictions \eqref{eq:m=3_e=1_transp} are satisfied for a positive proportion of specializations, there exists a transposition in $G$ by Corollary \ref{cor:hilbert_irred}. This concludes the proof of the lemma in the case that $f$ has a squarefree factor. 

Suppose now that $f = f_1^2$ for $f_1$ squarefree. This time, we find a prime $p\nmid c_d$ and an integer $\alpha_0$ such that $p^5 \ || \ f_1(\alpha_0)$ and $p \nmid f_1'(\alpha_0)$. This implies $p^{10} \ || \ f(\alpha_0)$, $p^5 \ || \ f'(\alpha_0)$, and $p \nmid f''(\alpha_0)$. Now consider an integral specialization $\a_0, \b_0, \balpha_0$ satisfying the following conditions:
\begin{align}
    v_p(a_0) &= 4 \\
    \nonumber v_p(a_1) &= 0 \\
    \nonumber v_p(a_i) & \geq 4 \text{ for all } 1 < i \leq \lfloor \frac{n}{3} \rfloor \\
    \nonumber v_p(b_j) &= 1 \text{ for all } 0 \leq j \leq (n-rd)/3 \\
    \nonumber v_p(\alpha_\ell) &= 0 \text{ for all } 0 < \ell \leq r.
\end{align}
As before, these conditions ensure that each coefficient of $F$ is divisible by $p^3$, except for the $t^3$ term. Inspecting the terms of $F$ in \eqref{eq:low_deg_terms_F} more carefully, we see that the constant, linear, and quadratic terms have valuations 12, (at least) 8, and 3, respectively, producing the Newton polygon below in Figure \ref{fig:n_odd_transp_v2}. 
\begin{figure}[ht]
    \centering
    \caption{$\NP_{\mathbb{Q}_{p}}(F_{\a_0,\b_0,\boldsymbol{\alpha_0}})$ with transposition}
    \label{fig:n_odd_transp_v2}
    \vspace{1ex}

    \begin{tikzpicture}[scale=0.75]
	   \draw[->, thick] (-0.2, 0) -- (10.5, 0);
	   \draw[->, thick] (0, -0.2) -- (0, 3.5);
	
	   \filldraw[black] (0, 3) circle (2pt) node[left] {$(0, 12)$};
	   \filldraw[black] (2, 0.75) circle (2pt) node[right] {$(2,3)$};
          \filldraw[black] (3, 0) circle (2pt) node[below] {$(3,0)$};
          \filldraw[black] (10, 0.75) circle (2pt) node[above] {$(n, 3)$};
	
	   \draw (0, 3) -- (2, 0.75);
	   \draw (3,0) -- (2, 0.75);
	   \draw (3,0) -- (10, 0.75);
    \end{tikzpicture}
    \end{figure}
Using the same argument via Lemma \ref{lem:NP_cycle} and Corollary \ref{cor:hilbert_irred}, the first segment of length 2 and slope $\frac{-9}{2}$ together with the hypotheses on $n$ grant us a transposition in $G$.
\end{proof}

Next, we find a cycle of prime length $q > \frac{n}{2}$ in $G$. For this, we recall a result of Breusch \cite{Breusch} which generalizes Bertrand's postulate for primes in arithmetic progressions (see also Moree's article \cite{Moree}, which includes a summary of similar theorems). In particular, for $n \geq 14$, there exist primes $q_1, q_2 \in (n/2, n)$ such that $q_i \equiv i \pmod{3}$. 

\begin{Lemma}\label{lem:m=3_q_cycle}
    Fix $m=3$ and an integral polynomial $f(x)$ of degree $d > 3$ with $3 \nmid d$ satisfying \eqref{eq:f_not_power}. If $n \geq \max(2d - 2, 14)$ and $3 \nmid n$, then the polynomial family $F(t)$ given in \eqref{eq:param_deg2} has a $q$-cycle in its Galois group for a prime $q > n/2$. 
\end{Lemma}

\begin{proof}
    By our hypotheses on $n$, we apply Breusch's result to ensure the existence of a prime $q \equiv n \pmod{3}$ in the interval $(n/2, n)$. 
    
    Consider now the integral specialization $\a_0, \b_0, \balpha_0$ satisfying the following conditions:
    \begin{align}
        v_p(a_i) &\geq 3 \text{ for all } i \neq \frac{n-q}{3} \\
        \nonumber v_p(a_{(n-q)/3}) &= 0 \\
        \nonumber v_p(b_j) &= 1 \text{ for all } 0 \leq j \leq (n-rd)/3 \\
        \nonumber v_p(\alpha_k) &= 0 \text{ for all } 0 \leq k \leq r.
    \end{align}
        This ensures that every coefficient of $F_{\a_0, \b_0, \balpha_0}$ is divisible by $p^3$ except for the $t^{n-q}$ term, resulting in the Newton polygon below in Figure \ref{fig:m=3_q_cycle}.

        \begin{figure}[ht]
    \centering
    \caption{$\NP_{\mathbb{Q}_{p}}(F_{\a_0,\b_0,\boldsymbol{\alpha_0}})$ with $q$-cycle}
    \label{fig:m=3_q_cycle}
    \vspace{1ex}

    \begin{tikzpicture}[scale=0.75]
	   \draw[->, thick] (-0.2, 0) -- (10.5, 0);
	   \draw[->, thick] (0, -0.2) -- (0, 2.5);
	
	   \filldraw[black] (0, 2) circle (2pt) node[left] {$(0, \geq 3)$};
	   \filldraw[black] (4,0) circle (2pt) node[below] {$(n-q,0)$};
          \filldraw[black] (10, 2) circle (2pt) node[right] {$(n, 3)$};
	
	   \draw[dashed] (0, 2) -- (4,0);
	   \draw (4,0) -- (10,2);
    \end{tikzpicture}
    \end{figure}
        Since the rightmost segment has length $q$ and slope $\frac{3}{q}$, Lemma \ref{lem:NP_cycle} together with Corollary \ref{cor:hilbert_irred} gives a $q$-cycle in $G$. 
\end{proof}

\begin{proof}[Proof of Proposition \ref{prop:m=3_irr_sn}]
    We have $F$ is irreducible by Proposition \ref{prop:param_irred_Sn}, so $G \subseteq S_n$ is transitive. Applying Lemmas \ref{lem:m=3_transp}, and \ref{lem:m=3_q_cycle}, we find that $G$ contains a transposition and a $q$-cycle for a prime $q > \frac{n}{2}$. We conclude the proof by an application of Lemma \ref{lem:genset} to find $G = S_n$.
\end{proof}

\begin{Example}\label{ex:running_example_V_galois}
    Let us return to our running concrete example (see Examples \ref{ex:running_example_I_S6point}, \ref{ex:running_example_II_S14point}, \ref{ex:running_example_III_higherdeg}, \ref{ex:running_example_IV_irr}),
    \[C \colon y^3 = x^4 + 1.\]
    Proposition \ref{prop:m=3_irr_sn} shows for $n \geq 14$ with $n \equiv 2,4 \pmod{6}$, the family $F(t)$ from Proposition \ref{prop:param} has Galois group $S_n$. This generalizes our observations in the previous examples from the special cases of $n=6,14$. Note that Proposition \ref{prop:m=3_irr_sn} is clearly not sharp: when $n=6$, for instance, we have seen that the family also has Galois group $S_6$.
\end{Example}

\begin{remark}
    It is not obvious to the authors how one might extend these proof techniques to other cases of interest, even in the trigonal case. Namely, it would be interesting to know that for the generic case of $3 \mid d$, when $n$ is a sufficiently large multiple of 3, the polynomial family $F(t)$ has Galois group $S_n$. In this case, the Newton polygon approach of Lemma \ref{lem:m=3_q_cycle} does not seem to easily produce a long cycle of prime length. 
    
    However, for general $m$ and sufficiently large $n$ in certain residue classes modulo $m$, the approach of Lemma \ref{lem:m=3_q_cycle} can be extended to find long cycles of prime length in $G$.

    In the event that $G$ contains a cycle of prime length in $(n/2, n-2)$, it follows that $G$ is a primitive subgroup of $S_n$ (i.e.\ it fixes no partition of the $n$ letters), and Jordan's theorem implies that $G$ contains the alternating group $A_n$. The approach of Lemma \ref{lem:m=3_transp} to show $G = S_n$ by producing a transposition in $G$ becomes unwieldy in general. An alternative approach is to show the discriminant of $F$ as a polynomial in $t$ is not squarefull, yielding a specialization whose Galois group contains a transposition; this approach was used in \cite{LOT}.
\end{remark}

\section{Lower bounds for $N_{n,C}(X)$}
\label{sec:lower_bounds}

In this section, we describe how to obtain the asymptotic lower bounds in Theorems \ref{thm:bound} and \ref{thm:m=3_bound}. We state this as a separate proposition to clarify that the counting argument is valid whenever $F$ is known to be irreducible over $\Q(\a, \b, \balpha)$ or have Galois group $S_n$.

\begin{Proposition}
\label{prop:bound}
    Fix an integer $m \geq 2$, a degree $d \geq m$ polynomial $f(x) \in \Z[x]$ satisfying \eqref{eq:f_not_power}, and an integer $n \geq n_0 = \max(d, \Frob(m,d) + 1)$. Suppose further that $F$ as given in Proposition \ref{prop:param} is irreducible. Then we have
    \[N_{n,C}(X) \gg X^{\delta_n},\]
    where $\delta_n$ is given by
    \begin{equation}\label{eq:c_n}
	\delta_n= \textstyle\frac{1}{m^2} + \frac{2n^2(m-m^2-dr+3)+n(km-k^2+4(m-m^2-dr)-dmr+d^2r^2)+2(km-k^2-dmr+d^2r^2)}{2m^2n^2(n-1)}.
\end{equation}    
    Here we take $r = \deg \gamma$ to be the minimal positive integer such that $n \equiv dr \pmod{m}$ as in Proposition \ref{prop:param}, and 
    \[k = \begin{cases}
        \min\{ k_1 \in \Z_{\geq 0} \mid \frac{n-d-k_1}{m} \in \Z\} & m \mid n, \\ 
        \min\{ k_2 \in \Z_{\geq 0} \mid \frac{n-k_2}{m} \in \Z\} & m \nmid n.       
    \end{cases}\]

    Moreover, if $\Gal(F/\Q(\a, \b, \balpha)) \simeq S_n$ then we have $N_{n,C}(X, S_n) \gg X^{\delta_n}$.
\end{Proposition}

The first statement of Theorem \ref{thm:bound} (resp.\ Theorem \ref{thm:m=3_bound}) follows from Proposition \ref{prop:bound} combined with Proposition \ref{prop:param_irred_Sn} (resp.\ Proposition \ref{prop:m=3_irr_sn}). The second statement, an improvement to the exponent when $n$ is sufficiently large, follows from the discussion in Section \ref{sec:improve}, namely Corollary \ref{cor:bound_improvements}.

\subsection{Coefficient bounds}

In this section, we construct a family of polynomials $P_{f,n}(Y)$ arising from certain specializations of \eqref{eq:param} in \Cref{sec:parametrization} satisfying bounds on the coefficients of $g(t)$ and $h(t)$ in $F(t)=g(t)^m-h(t)^mf(\gamma(t))$. These bounds will be useful for counting multiplicities of fields generated by these family of polynomials, and help control the discriminant via the following well known consequence of the homogeneity of the discriminant.

\begin{Lemma}
\label{lem:roots_coeff_bound}
	Let $F(x) = \sum_{i=0}^n d_ix^{i} \in \C[x]$ be a polynomial of degree $n$ and suppose $Y > 0$. If for all $i < n$ we have $\left| d_i \right| \leq Y^{n-i}$, then $\left|\Disc(F)\right| = O(Y^{n(n-1)})$, where the implied constant depends only on $n$ and the leading term $d_n$.
\end{Lemma}

Before defining $P_{f,n}(Y)$ in general as promised, we give two examples.

\begin{Example}\label{ex:running_example_VI_countsetup_6}
    Let us return to our running concrete example (see Examples \ref{ex:running_example_I_S6point}, \ref{ex:running_example_II_S14point}, \ref{ex:running_example_III_higherdeg}, \ref{ex:running_example_IV_irr}, \ref{ex:running_example_V_galois}),
    \[C \colon y^3 = x^4 + 1.\]    
    Fix $n=6$. We know from Example \ref{ex:running_example_I_S6point} and Hilbert's irreducibility theorem that the partial specialization $F_{a_{2}=1,\gamma_0 = t}$ is irreducible. Consider for $Y > 0$ the set of partial specializations
    \[P_{x^4+1, 6}(Y) = \left\{F = b_0^3(t^4 + 1) - (t^2 + a_1t + a_0)^3 : \left|a_1\right| \leq Y,\ \left|a_0\right| \leq Y^2,\ \left|b_0\right| \leq Y^{2/3} \right\}.\]
    All $F \in P_{x^4+1, 6}(Y)$ are sextic polynomials with leading coefficient $-1$. Writing $F= \sum d_it^i$, the bounds on the coefficients force $\left| d_i \right| \ll Y^{6-i}$. Applying Lemma \ref{lem:roots_coeff_bound}, we have $\left|\Disc(F)\right| \ll Y^{30}$ for $F \in P_{x^4+1, 6}(Y)$.
\end{Example}

\begin{Example}\label{ex:running_example_VII_countsetup_14}
    Continuing with $C$ from our previous example, but setting $n=14$, we build $P_{x^4+1, 14}(Y)$ similarly. Thanks to Example \ref{ex:running_example_II_S14point}, we have that the partial specialization $F_{b_2=1,\gamma_0 = t^2}$ is irreducible. Consider for $Y > 0$ the set of partial specializations
    \begin{align*}
        P_{f,n}(Y) = \Big\{F = (t^2 + b_1t + b_0)^3(t^8 + 1) - (a_4t^4 + a_3t^3 &+ a_2t^2 + a_1t + a_0)^3 : \\
        &\left|a_{4-i}\right| \leq Y^{2/3 + i},\ \left|b_1\right| \leq Y,\ \left|b_0\right| \leq Y^2 \Big\}.
    \end{align*}
    Again writing $F= \sum d_it^i$, the bounds on the coefficients above force $\left| d_i \right| \ll Y^{14-i}$, so Lemma \ref{lem:roots_coeff_bound} gives $\left|\Disc(F)\right| \ll Y^{182}$ for $F \in P_{x^4+1, 14}(Y)$.
\end{Example}  

The construction of $P_{f,n}(Y)$ in general follows the approach of Examples \ref{ex:running_example_VI_countsetup_6} and \ref{ex:running_example_VII_countsetup_14}: choose $\gamma_0$ and the leading coefficients of $g$ and/or $h$ so that the partial specialization of $F$ is irreducible, then restrict the coefficients $a_i, b_j$ so that we can apply Lemma \ref{lem:roots_coeff_bound} to bound the discriminant. For the remainder of this section, we work within the hypotheses of Proposition \ref{prop:bound}, namely that the polynomial family $F \in \Q(\a, \b, \balpha)[t]$ is irreducible. When $F$ is further known to have Galois group $S_n$, the argument used to bound $N_{n,C}(X)$ produces the same bound on $N_{n,C}(X, S_n)$.

First, we argue that we can find appropriate partial specializations.

\begin{Lemma}\label{lem:partial_specialization}
    Assume the same hypotheses as Proposition \ref{prop:bound}. Then there exists $\balpha_0 \in \Z^{r+1}$ such that the partial specialization $F_{\balpha_0} \in \Q(\a, \b)[t]$ is irreducible and $f(\gamma_0(t))$ is also $m$-th power free.

    Moreover, if $\Gal(F/\Q(\a, \b, \balpha)) \simeq S_n$ then $\balpha_0$ may be chosen such that the partial specialization also has full Galois group, $\Gal(F_{\balpha_0}/\Q(\a, \b)) \simeq S_n$.
\end{Lemma}

\begin{proof}
    $F$ is irreducible over $\Q(\a, \b, \balpha)$ by hypothesis, so Hilbert's irreducibility theorem (Lemma \ref{thm:hilbert_irred} for arbitrary base field) implies that for almost all choices of $\balpha_0$, $F_{\balpha_0}$ is irreducible over $\Q(\a,\b)$. If $\Gal(F/\Q(a,b,\balpha)) \simeq S_n$ then for almost all $\balpha_0$ we have $\Gal(F_{\balpha_0}/\Q(\a,\b)) \simeq S_n$. 

    Let us now examine more closely when $f(\gamma_0(t))$ is $m$-th power free. Write $f$ as a product of irreducible factors $f = \prod f_i^{e_i}$ and recall that $f$ itself is $m$-th power free. Setting $f_{\mathrm{rad}} = \prod f_i$, it is enough to check that we can choose $\gamma_0$ such that $f_{\mathrm{rad}}(\gamma_0(t))$ is \textit{squarefree}.

    Now we may use the discriminant $\Disc\left(f_{\mathrm{rad}}(\gamma(t))\right)$, viewed as a polynomial function in the variables $\balpha$; $f_{\mathrm{rad}}(\gamma(t))$ has a multiple root wherever this polynomial vanishes, which is a Zariski closed condition on the affine space $\mathbb{A}^{r+1}$ from which we are choosing $\balpha_0$. Hence the space of $\balpha_0$ giving rise to $f_{\mathrm{rad}}(\gamma_0(t))$ that are squarefree --- and thus to $f(\gamma_0(t))$ that are $m$-th power free --- is Zariski dense. In particular, some such $\balpha_0$ satisfies both the condition that $F_{\balpha_0}$ is irreducible and that $f(\gamma_0(t))$ is $m$-th power free.    
\end{proof}

Moving forward, we fix some $\gamma_0(t)$ such that $F_{\balpha_0}$ is irreducible (with Galois group $S_n$ if appropriate) and such that $f(\gamma_0(t))$ {is} $m$-th power free, by Lemma \ref{lem:partial_specialization}. 

Let $Y > 0$ be a real number. When $m \mid n$ we write
\begin{align}\label{eq:restrictions_m|n}
	g(t) &= a_{n/m}t^{n/m}+a_{n/m-1}t^{n/m-1}+ \dots + a_0,\\
	\nonumber h(t)&= b_{(n-d-k_1)/m}t^{(n-d-k_1)/m}+b_{(n-d-k_1)/m-1}t^{(n-d-k_1)/m-1}+ \dots +b_0,
\end{align}  
Here $k_1$ is the minimal nonnegative integer such that $(n-d-k_1)/m$ is an integer. This realizes the degrees in \eqref{eq:param_deg1}.

In the case where $m\nmid n$, we choose  $r$, the degree of $\gamma(t)$,  to be the minimal positive integer for which $n \equiv dr \pmod{m}$. As above, we have
\begin{align}\label{eq:restrictions_mnmidn}
	g(t) &= a_{(n-k_2)/m}t^{(n-k_2)/m}+a_{(n-k_2)/m-1}t^{(n-k_2)/m-1}+ \dots + a_0,\\
	\nonumber h(t) &= b_{(n-dr)/m}t^{(n-dr)/m}+ b_{(n-dr)/m-1}t^{(n-dr)/m-1}+ \dots + b_0, 
\end{align}
  Here $k_2$ is the minimal positive integer such that $(n-k_2)/m$ is an integer so this realizes the degrees in \eqref{eq:param_deg2}.

\begin{Definition}[$P_{f,n}(Y)$]
\label{def:PfnY}
    For $Y > 0$, we define 
    \[P_{f,n}(Y) = \left\{ F(t) = h(t)^m f(\gamma_0(t)) - g(t)^m \right\}\]
    where $g,h$ satisfy the following conditions:
    \begin{itemize}
        \item if $m \mid n$, $g,h$ are given by \eqref{eq:restrictions_m|n} with 
        \begin{itemize}
            \item $\left|a_{n/m-i}\right| \leq Y^i$ for $i > 0$, 
            \item $\left|b_{(n-d-k_1)/m - j}\right| \leq Y^{k_1/m + j}$ for $j > 0$
            \item $a_{n/m}$ is a fixed integer so that $F_{a_{n/m}, \balpha_0}$ is irreducible,
            \item if $k_1 = 0$ then $b_{(n-d)/m}$ is fixed so that $F_{a_{n/m}, b_{(n-d)/m}, \balpha_0}$ is irreducible, and if $k_1 \neq 0$ then $\left|b_{(n-d-k_1)/m} \right| \leq Y^{k_1/m}$;
        \end{itemize}
        \item if $m \nmid n$, $g,h$ are given by \eqref{eq:restrictions_mnmidn} with 
        \begin{itemize}
            \item $\left|a_{(n-k_2)/m - i}\right| \leq Y^{k_2/m + i}$ for $i \geq 0$, 
            \item $\left|b_{(n-dr)/m - j}\right| \leq Y^{j}$ for $j > 0$
            \item $b_{(n-dr)/m}$ is a fixed integer so that $F_{b_{(n-dr)/m}, \balpha_0}$ is irreducible.
        \end{itemize}
    \end{itemize}
\end{Definition}

Note that the existence of $a_{n/m}$ (and $b_{(n-d)/m}$ if $k_1 = 0$) or $b_{(n-dr)/m}$ so the appropriate partial specialization is irreducible is ensured by Theorem \ref{thm:hilbert_irred}.

The polynomials $F \in P_{f,n}(Y)$ have degree $n$. Writing $F(t) = d_nt^n + \ldots + d_0$, we see that the restrictions on the coefficients in Definition \ref{def:PfnY} ensure that $\left|d_i\right| \ll Y^{n-i}$. Applying Lemma \ref{lem:roots_coeff_bound}, we see that for all $F \in P_{f,n}(Y)$, we have $|\text{Disc}(F)| \leq BY^{n(n-1)}$ for a constant $B$ depending on $m$, $f$, and $n$. 

\subsection{Bounding multiplicities}
We count the number of fields arising from specializations in  \eqref{eq:param} by counting the number of polynomials in $P_{f,n}(Y)$ and adjusting for two possible sources of multiplicity. The first potential source of multiplicity is the case where two different $g(t)$, $h(t)$ give rise to the same element $F(t)$ in $P_{f,n}(Y)$. The second potential source of multiplicity is that multiple elements $F(t)$ in $P_{f,n}(Y)$ produce isomorphic number fields. The first potential source of multiplicity is dealt with by the following lemma, building on the strategy in \cite[Lemma 7.4]{LOT}.

\begin{Lemma}
\label{lem:poly_mult}
	Let $F(t) \in \Z[t]$ be a polynomial of degree $n$. The number of ways to choose nonzero polynomials $g(t), h(t) \in \Z[t]$ of some fixed degrees $\deg g \leq \frac{n}{m}$ and $\deg h < \frac{n}{m}$ with one of the leading coefficients of $g$ or $h$ fixed, such that $F(t) = g(t)^m - f(\gamma_0(t))h(t)^m$ is $O_{m,n}(1)$.
\end{Lemma}

\begin{proof}
	Note that we assumed $f(x)$ is $m$-th power free in \eqref{eq:f_not_power}. We then chose $\gamma_0(t)$ as in Lemma \ref{lem:partial_specialization} such that $f(\gamma_0(t))$ is also $m$-th power free. The coordinate ring $R = \C[t,y]/(y^m-f(\gamma_0(t)))$ is a Noetherian domain of Krull dimension one, thus its integral closure $\wt{R}$ is a Dedekind domain. This implies that in $\wt{R}$, the ideal $(F)$ factors uniquely into a product of finitely many primes, of the form $(t-t_0, y-y_0)$ satisfying both $y_0 = f(\gamma(t_0))$ and $F(t_0) = 0$. There are $mn$ such solutions, counted with multiplicity, so we have at most $mn$ prime factors of $(F)$.
	
	 As in the proof of \cite[Lemma 7.4]{LOT}, we observe that given any such $g,h$ there is a factorization 
	\[F = g^m - f(\gamma_0)h^m = \prod_{i=0}^{m-1} \left(g-\zeta^if(\gamma_0)^{1/m} h\right),\]
	where $\zeta$ is a primitive $m$-th root of unity. The ideal $(g-f(\gamma_0)^{1/m}h)$ divides $(F)$ so there are at most finitely many possibilities for its prime factorization. Thus there are at most finitely many choices for the ideal $(g-f(\gamma_0)^{1/m}h)$. It remains to show that if $g$ and $h$ satisfy the hypotheses of the lemma, this ideal determines them precisely.
	
	Suppose we have $g', h'$ satisfying the hypotheses with $(g-f(\gamma_0)^{1/m}h) = (g'-f(\gamma_0)^{1/m}h')$. Then for some unit $u\in ({\wt{R}})^\times$, we have $g-f(\gamma_0)^{1/m}h = u(g'-f(\gamma_0)^{1/m}h')$. This unit $u$ necessarily satisfies a minimal monic polynomial
	\begin{equation}
	\label{eq:min_poly_for_unit}
		u^k + v_{k-1}u^{k-1} + \cdots + v_1u + v_0 = 0,
	\end{equation}
	where $v_i \in R$ and $v_0 \in R^\times$. Multiplying by $g' - f(\gamma_0)^{1/m}h'$, this becomes
	\begin{align*}
		0 &= \left(g' - f(\gamma_0)^{1/m}h'\right) \left( u^k + v_{k-1}u^{k-1} + \cdots + v_1u + v_0 \right) \\ 
		&= (g - f(\gamma_0)^{1/m}h)u^{k-1} + (g - f(\gamma_0)^{1/m}h)v_{k-1}u^{k-1} + \cdots + (g - f(\gamma_0)^{1/m}h)v_1 + v_0\left(g' - f(\gamma_0)^{1/m}h'\right).
	\end{align*}
	If $k > 1$ then this contradicts minimality of \eqref{eq:min_poly_for_unit}, so we must have $k=1$, in which case we have $u \in R^\times$.
	
	With this in hand, we may write $u = u(t)$ as
	\[u(t) = \sum_{i=0}^{m-1} f(\gamma_0(t))^{i/m} u_i(t)\]
	with $u_i(t) \in \C[t]$. The relation $u(g-f(\gamma_0)^{1/m}h) = g'-f(\gamma_0)^{1/m}h'$ implies that
    \begin{align}\label{eq:u_system}
		u_0g - u_{m-1}f(\gamma_0)h &= g' \\
		\nonumber u_1g - u_0h &= h' \\
		\nonumber u_ig - u_{i-1}h &= 0 \quad \text{for all }2 \leq i \leq m-1
	\end{align}
	as polynomials in $\C[t]$. Tracing through \eqref{eq:u_system}, we determine that
	\begin{equation}\label{eq:u_m-1}
		u_{m-1}F = u_{m-1}(g^m - h^mf(\gamma_0)) = g'h^{m-1} + gh^{m-2}h'.
	\end{equation}
	If $u_{m-1} \neq 0$, the left hand side has degree $\deg u_{m-1} + n$, while the right hand side has degree at most $\deg g + (m-1) \deg h < n$, producing a contradiction. Therefore, we have $u_{m-1} = 0$, and tracing through the relations \eqref{eq:u_system} again, this implies $u_i = 0$ for all $1 \leq i \leq m-1$, i.e.\ $u(t) = u_0(t)$.
	
	Finally, we observe that since the degrees of $g$ and $h$ are fixed, $u = u_0$ must be a constant. Moreover, since we also require the leading coefficients of either $g,g'$ or $h,h'$ to be fixed, we must have $u=1$. Therefore, the ideal $(g-f(\gamma_0)^{1/m}h)$ can come from at most one $g,h$ satisfying the hypotheses.
\end{proof}

When $m \mid n$, the restrictions imposed in \eqref{eq:restrictions_m|n} and the definition of $P_{f,n}(Y)$ fix the degrees of $g$ and $h$ and the leading coefficient of $g$  such that the hypotheses of Lemma \ref{lem:poly_mult} are satisfied. Thus each choice of $g(t)$ and $h(t)$ coincides with at most finitely many others. The same is true for the $m \nmid n$ case. Thus we can give a count for the number of $F(t)$ in $P_{f,n}(Y)$ based on the number of choices for $g(t)$ and $h(t)$. More precisely, $\# P_{f,n}(Y) \asymp Y^c$ where $c$ is defined as follows.

In the case where $m\mid n$, we have
\begin{equation}\label{eq:c_m|n} c= \sum\limits_{i=1}^{n/m} i + \sum\limits_{j=0}^{(n-d-k_1)/m} \left(j+\dfrac{k_1}{m}\right)= \dfrac{1}{m^2} 
\left( n^2+ n(m-d)+ \dfrac{d^2+(k_1-d)m-k_1^2}{2} \right). \end{equation}

\noindent In the case where $m\nmid n$, we have
\begin{equation}\label{eq:c_m_nmid_n} c=\sum\limits_{i=0}^{(n-k_2)/m} \left(\dfrac{k_2}{m}+i\right) + \sum\limits_{j=1}^{(n-rd)/m} j = \dfrac{1}{m^2} \left(n^2+n(m-dr)+\dfrac{d^2r^2+(k_2-dr)m-k_2^2}{2} \right).\end{equation}
Let $P_{f,n}(Y, \irr)$ denote the subset of $P_{f,n}(Y)$ consisting of irreducible polynomials (similarly, let $P_{f,n}(Y, S_n)$ denote the subset of irreducible polynomials with Galois group $S_n$ over $\Q$). Since we have assumed $F$ is irreducible and chosen $a_{n/m}$ or $b_{(n-dr)/m}$ appropriately, Lemma \ref{lem:partial_specialization} implies that $\# P_{f,n}(Y, \irr) \asymp Y^c$ (and similarly $\# P_{f,n}(Y, S_n) \asymp Y^c$ if $F$ has symmetric Galois group).

We now turn to the other source of multiplicity: distinct polynomials cutting out the same number field.

\begin{Example}\label{ex:running_example_VIII_multiplicity}
    Let us return to our running concrete example (see Examples \ref{ex:running_example_I_S6point}, \ref{ex:running_example_II_S14point}, \ref{ex:running_example_III_higherdeg}, \ref{ex:running_example_IV_irr}, \ref{ex:running_example_V_galois}, \ref{ex:running_example_VI_countsetup_6}, \ref{ex:running_example_VII_countsetup_14}),
    \[C \colon y^3 = x^4 + 1.\]    
    Recall from Example \ref{ex:running_example_I_S6point} that for $n=6$, taking $g=t^2-t$, $h=1$, $\gamma = t$, we obtained $F = -t^6 + 3t^5 - 2t^4 + t^3 + 1$. Having chosen $a_2 = 1$, $\gamma_0 = t$ in our construction in Example \ref{ex:running_example_VI_countsetup_6}, we have $F \in P_{x^4+1,6}(Y)$ for all $Y \geq 1$. Set $K = \Q[t]/F(t)$.

    If we instead take $g'=t^2+t$ (leaving $h, \gamma$ unchanged), we have 
    \[F' = t^4+1 - (t^2-t)^3 = -t^6 - 3t^5 - 2t^4 - t^3 + 1 \in P_{x^4+1,6}(Y)\] for all $Y \geq 1$. A computation reveals $K \simeq \Q[t]/F'(t)$, illustrating that two polynomials in $P_{x^4+1,6}(Y)$ can cut out the same field.
\end{Example}

To address this source of potential multiplicity (that there may be multiple elements of $F(t)$ that produce isomorphic number fields), we build on a strategy of Ellenberg and Venkatesh \cite{ellenbergvenkatesh} for counting number fields and appeal to the multiplicity counts of Lemke Oliver and Thorne \cite{LOT}. See also \cite{keyes} for a more detailed discussion.

Our restrictions on the sizes of $|a_i|,\ |b_j|$ in the definition of $P_{f,n}(Y)$ ensure that the coefficients of $F$ satisfy $|d_{n-i}| \leq AY^i$ for some constant $A$. In particular the leading terms are fixed, hence we may divide by them to obtain monic polynomials.

\begin{Definition}[$S(Y)$]
    Let $Y > 0$ and define
    \[S(Y) := \set{ F = t^n + d_{n-1}'t^{n-1} + ... + d_0' \in (1/w)\Z[t] : F \text{ irreducible and } \left|d_{n-i}' \right| \ll_{n,f} Y^i }\]
\end{Definition}

Note that by this construction, elements of $P_{f,n}(Y, \irr)$ (and $P_{f,n}(Y, S_n)$) are in bijection with a subset of $S(Y)$, provided we choose the implied constant appropriately. 

We define the multiplicity of a number field $K$ of degree $n$ in $S(Y)$ to be the number of polynomials in $S(Y)$ that cut out the field $K$,
\[M_K(Y) := \# \set{F \in S(Y) \mid \Q[t]/F(t) \simeq K}.\]
We state here several useful bounds related to this multiplicity $M_K(Y)$ that we will need. The following is due to Lemke Oliver and Thorne. 

\begin{Lemma}[Lemke Oliver--Thorne {\cite[Proposition 7.5]{LOT}}]
\label{lem:LOT_mult}
	We have \[M_K(Y) \ll \max\Big(Y^n\left|\Disc(K)\right|^{-1/2}, Y^{n/2}\Big).\]
\end{Lemma}

Lemma \ref{lem:LOT_mult}, together with the following upper bound for $N_{n}(X)$ due to Schmidt, is used in \cite{keyes} to give a bound for the sum of multiplicities of fields with discriminant bounded by $T$.

\begin{Theorem}[Schmidt, \cite{schmidt}]
\label{prop:counting_fields_general}
	For $n \geq 3$, we have
	\begin{equation}\label{eq:schmidt_bound}
		N_n(X) \ll X^{\frac{n + 2}{4}}.
	\end{equation}
\end{Theorem}

\begin{Lemma}[Keyes, {\cite[Lemma 5.4]{keyes}}]
\label{lem:sum_small_disc}
	Let $T \leq Y^n$. Then \[\sum_{\left|\Disc(K)\right| \leq T} M_K(Y) \ll Y^nT^{n/4},\] where the sum runs over all degree $n$ number fields $K$ such that $\left|\Disc(K)\right| \leq T$.
\end{Lemma}

\begin{remark}
        Theorem \ref{prop:counting_fields_general} has been superseded for $n \geq 6$, but we use it anyway for now due to its simple form. We defer further discussion of improvements to Lemma \ref{lem:sum_small_disc} until Section \ref{sec:improve}, where we discuss how better upper bounds for $N_{n}(X)$ improve our lower bounds on $N_{n,C}(X)$ for $n$ sufficiently large.
\end{remark}

\subsection{Bounding $N_{n,C}(X)$}

We now have all the tools to complete the proof of Proposition \ref{prop:bound}. We first walk through the proof in an example.

\begin{Example}\label{ex:running_example_IX_result}
    Let us return to our running concrete example (see Examples \ref{ex:running_example_I_S6point}, \ref{ex:running_example_II_S14point}, \ref{ex:running_example_III_higherdeg}, \ref{ex:running_example_IV_irr}, \ref{ex:running_example_V_galois}, \ref{ex:running_example_VI_countsetup_6}, \ref{ex:running_example_VII_countsetup_14}, \ref{ex:running_example_VIII_multiplicity}),
    \[C \colon y^3 = x^4 + 1.\]    
    Consider the case of $n=14$. In Example \ref{ex:running_example_VII_countsetup_14} we defined $P_{x^4+1,14}(Y)$, a set of partial specializations of the polynomial family from Proposition \ref{prop:param}. By Lemma \ref{lem:poly_mult}, at most $O(1)$ choices of $g,h$ result in coinciding $F$, so we have 
    \[\#P_{x^4+1,14}(Y) \asymp Y^{49/3},\]
    Lemma \ref{thm:hilbert_irred} allows us to further conclude $\#P_{x^4+1,14}(Y, \irr) \asymp \#P_{x^4+1,14}(Y, S_{14}) \asymp Y^{49/3}$. Since $P_{x^4+1,14}(Y)$ already consists of monic polynomials, we have $S(Y) = P_{x^4+1,14}(Y, \irr)$ and $S(Y, S_{14}) = P_{x^4+1,14}(Y, S_{14})$.

    The upper bound on $M_K(Y)$ of Lemma \ref{lem:LOT_mult} decreases in $\left|\Disc(K)\right|$, so when accounting for multiplicity, we would like the contribution from $F \in S(Y)$ with small discriminant to be small. To make this precise, we set $T = \kappa Y^{2/3}$ and apply Lemma \ref{lem:sum_small_disc} to see that
    \[\sum_{\left|\Disc(K) \right| \leq Y^{2/3}} M_K(Y) \ll \kappa^{14/4} Y^{49/3}.\]
    After choosing $\kappa$ appropriately, we have that at most a positive proportion of $F \in S(Y)$ yield number fields with discriminant at most $T$. Thus we focus instead on those with discriminant between $T$ and a constant times $Y^{13\cdot 14} = Y^{182}$.

    For these fields $K$, Lemma \ref{lem:LOT_mult} shows
    \[M_K(Y) \ll Y^{41/3}.\]
    We can bound the number of such fields $K$ simply by taking our asymptotic count of $F \in S(Y)$ and dividing by $Y^{41/3}$:
    \[N_{14,C}(BY^{182}) \gg \sum_{\substack{K = \Q[t]/F(t) \text{ for } F \in S(Y) \\ T \leq \left| \Disc(K) \right| \leq BY^{182}}} 1 \gg \frac{\#S(Y)}{Y^{41/3}} \gg \frac{Y^{49/3}}{Y^{41/3}} = Y^{8/3}.\]
    Setting $X = BY^{182}$ yields $N_{14,C}(X) \gg X^{4/273}$, which is seen to agree with the value of $\delta_{14}$ of Theorem \ref{thm:bound}, given in \eqref{eq:final_bound}, with $m=3$, $d=4$, $n=14$, $k=k_2=2$, $r=2$.
\end{Example}

\begin{proof}[Proof of Proposition \ref{prop:bound}]
By our construction, for any $F \in P_{f,n}(Y, \irr)$ and any root $\alpha$ of $F$, we have $(\alpha, \frac{g(\alpha)}{h(\alpha)}) \in C(K)$ where $K = \Q(\alpha)$ is a field of degree $n$. Recall also that we have $\left|\Disc(K)\right| \leq BY^{n(n-1)}$ for a constant $B$. 

First we will show that fields of low discriminant here are negligible in their contributions to $N_{n,C}(X)$. Using Lemma \ref{lem:sum_small_disc}, we choose $T= \kappa Y^{\frac{1}{m^2}(4n-4(dr+(m-1)m)+(2(dr-k)(dr+k-m))/n)}$ so that

\begin{equation}
\label{eq:small_disc}
    \sum_{\left|\Disc(K)\right| \leq T} M_K(Y) \ll \kappa^{n/4}Y^c,
\end{equation} 
and we recall that 
\begin{equation}
\label{eq:poly_asymp}
		\#P_{f,n}(Y, \irr) \asymp Y^c
\end{equation} 
where $c$ is given either by \eqref{eq:c_m|n} or \eqref{eq:c_m_nmid_n}. We choose $\kappa$ to be sufficiently small so that the quantity in \eqref{eq:small_disc} is at most $\# P_{f,n}(Y, \irr)/2$. Thus our parameterization produces negligibly many fields of discriminant at most $T$. Since the bound in Lemma \ref{lem:LOT_mult} is decreasing with respect to $|\Disc(K)|$, thus we have $M_K(Y) \ll T^{-1/2}Y^n$ for all $K$ of discriminant $T <\left|\Disc(K)\right| \leq BY^{n(n-1)}$. We obtain an asymptotic lower bound for $N_{n,C}(BY^{n(n-1)})$ by dividing $\#P_{f,n}(Y, \irr)$ by this worst case multiplicity.
\begin{align}\label{eq:almost_done}
	\nonumber N_{n,C}(BY^{n(n-1)}) &\gg Y^{c-n}T^{1/2}\\
	&=  \scalebox{0.93}{$\displaystyle Y^{\frac{1}{m^2}\left(n^2+n(2+m-m^2-dr)+(-k^2+4m+km-4m^2-4dr-dmr+d^2r^2)/2+(d^2r^2-k^2+km-dmr)/n\right)}$}.
\end{align}

To obtain the exponent $\delta_n$ in \eqref{eq:c_n}, we replace $Y$ in \eqref{eq:almost_done} by $\left(X/B\right)^{1/n(n-1)}$. This produces
\begin{equation*}
	\delta_n= \textstyle\frac{1}{m^2} + \frac{2n^2(m-m^2-dr+3)+n(km-k^2+4(m-m^2-dr)-dmr+d^2r^2)+2(km-k^2-dmr+d^2r^2)}{2m^2n^2(n-1)}
\end{equation*}
and thus $N_{n,C}(X) \gg X^{\delta_n}$, as desired. 

In the case $\Gal(F/\Q(\a, \b, \balpha_0)) \simeq S_n$, we repeat the argument by replacing $P_{f,n}(Y, \irr)$ by $P_{f,n}(Y, S_n)$ to obtain $N_{n,C}(X, S_n) \gg X^{\delta_n}$.
\end{proof}

\subsection{Improvements for $n$ sufficiently large}
\label{sec:improve}

As in \cite[Section 5.4]{keyes}, we can improve on our lower bound when $n$ is sufficiently large by employing better known upper bounds for $N_n(X)$. The idea is to show that if the upper bound for $N_n(X)$ is good enough, then the best case scenario of Lemma \ref{lem:LOT_mult} applies, and we can assume $M_K(Y) \ll Y^{n/2}$. Thus
\[N_n(Y^{n(n-1)}) \gg Y^{c - \frac{n}{2}}\]
where $c$ is given in \eqref{eq:c_m|n} or \eqref{eq:c_m_nmid_n}, as appropriate. It remains to compute this exponent and determine when the improved upper bounds for $N_n(X)$ take effect.

Assume we have an upper bound of the form
\[(*) \quad N_n(X) \ll X^{\varepsilon(n,m,d)},\]
where $\varepsilon(n,m,d) \geq 1$ is a constant depending on $n$ and the $m,d$ values for our curve $C$. We will use a modification of (the proof of) Lemma 5.5 which is somewhat more flexible. 
\begin{Lemma}\label{lem:sum_small_disc_v2}
	Let $T \leq Y^n$. Assume $(*)$ for some constant $\varepsilon(n,m,d)$. Then
		\[\sum_{|\Disc K| \leq T} M_K(Y) \ll Y^{n}T^{\varepsilon(n,m,d)-1/2} + \frac{Y^nT^{\varepsilon(n,m,d) - \frac{1}{2}}}{2\varepsilon(n,m,d) - 1}.\]
	In particular, when we take $T = Y^n$ we have
	\[\sum_{|\Disc K| \leq Y^n} M_K(Y) \ll Y^{\frac{n}{2} + n\varepsilon(n,m,d)}.\]
\end{Lemma}

\begin{proof}
Write \[M(Y)(t) = \max \left\{ M_K(Y) : |\Disc K| = t \right\}\]
	for the maximal multiplicity of a number field with discriminant $t$. Note that the bound in Lemma \ref{lem:LOT_mult} depends only on the discriminant so we have $M(Y)(t) \ll \max\left(Y^n t^{-1/2}, Y^{n/2}\right).$ We set up a Riemann-Stieljes integral as in \cite[Lemma 5.4]{keyes}. 
	\begin{align*}
		\sum_{|\Disc K| \leq T} M_K(Y) &\leq \int_{1^-}^T M(Y)(t) dN_n(t)\\
		&\ll  \int_{1^-}^T Y^nt^{-\frac{1}{2}}dN_n(t)\\
		&=Y^nT^{-\frac{1}{2}}N_n(T) + \frac{Y^n}{2}\int_{1^-}^T t^{-\frac{3}{2}} N_n(t) dt.
	\end{align*}
	Substituting $(*)$ into the last line above gives the first statement of the lemma.
\end{proof}

Note that Lemma \ref{lem:sum_small_disc} follows from this by taking $\varepsilon(n,m,d) = \frac{n+2}{4}$, the Schmidt bound \cite{schmidt}. However, this is not good enough for $Y^{\frac{n}{2} + n\varepsilon(n,m,d)}$ to be $o(Y^c)$. For this we need
\[(**) \quad \varepsilon(n,m,d) < \frac{c}{n} - \frac{1}{2}.\]
Using the best known upper bounds we can find when $(**)$ is satisfied for a given $C$ and $n$.

\begin{Theorem}[Lemke Oliver--Thorne, {\cite[Theorem 1.1]{LOT_upper}}]
\label{thm:LOT_upper}
	For $n \geq 6$ we have 
	\[N_n(X) \ll X^{1.564 (\log n)^2}.\]
\end{Theorem}

This is sufficient to give \eqref{eq:final_improved} in Theorem \ref{thm:bound} and the second statement of Theorem \ref{thm:m=3_bound}, both of which follow from the corollary below.

\begin{Corollary}\label{cor:bound_improvements}
    Fix an integer $m \geq 2$, a degree $d \geq m$ polynomial $f(x) \in \Z[x]$ satisfying \eqref{eq:f_not_power}. Then for all sufficiently large multiples $n$ of $\gcd(m,d)$, whenever $F$ as given in Proposition \ref{prop:param} is irreducible, we have
	\[N_n(X) \gg X^{\delta_n'},\]
	where \[\delta_n' ={\frac{1}{m^2} \left(1 + \frac{(2m - 2dr + 1)n + d^2r^2 - mdr + mk - k^2}{2n(n-1)}\right)}.\]
	This is precisely \eqref{eq:final_improved}.

    Moreover, whenever $n$ is sufficiently large and $\Gal(F/\Q(\a, \b, \balpha)) = S_n$, we have $N_{n,C}(X,S_n) \gg X^{\delta_n'}$.
\end{Corollary}

\begin{proof}
	Fix a choice of $C$, so $m$ and $d$ are fixed. Assume $n \geq 6$ and set $\varepsilon(n,m,d) = 1.564(\log n)^2$, so Theorem \ref{thm:LOT_upper} ensures $(*)$ is satisfied. Recalling $c$ from \eqref{eq:c_m|n} or \eqref{eq:c_m_nmid_n} we see that in either case, $\frac{c}{n} - \frac{1}{2}$ grows linearly with $n$, as $k_1, k_2,$ and/or $r$ are bounded, depending on $m,d$. Clearly $(\log n)^2$ grows more slowly with $n$, so for $n$ sufficiently large $(**)$ is satisfied.
	
	As noted above, Lemma \ref{lem:sum_small_disc_v2} together with $(*), (**)$ implies that 
	\[\sum_{|\Disc K| \leq Y^n} M_K(Y) = o(Y^c).\]
	Thus the contribution of fields with discriminant up to $Y^n$ to $\#P_{f,n}(Y,\irr)$ is negligible. For fields $K$ with $Y^n <|\Disc K| \leq Y^{n(n-1)}$ we have $M_K(Y) \ll Y^{n/2}$ by Lemma \ref{lem:LOT_mult}. Hence, we have
	\[N_{n,C}(Y^{n(n-1)}) \gg \#P_{f,n}(Y,\irr) Y^{-\frac{n}{2}} \gg Y^{c-\frac{n}{2}}.\]
	To get $\delta_n'$ we set $Y = X^{\frac{1}{n(n-1)}}$ and take $\delta_n' = \frac{c - n/2}{n(n-1)}$, which we can compute explicitly to obtain the stated value.

    In the case $\Gal(F/\Q(\a,\b,\balpha)) \simeq S_n$, we may replace $P_{f,n}(Y, \irr)$ by $P_{f,n}(Y, S_n)$ to obtain the lower bound $N_{n,C}(X, S_n) \gg X^{\delta_n'}$.
\end{proof}

The question remains to find when the improved asymptotic lower bound above takes effect. That is, to determine when $(*)$ and $(**)$ are both satisfied. To do this, we make use of a more flexible version of Theorem \ref{thm:LOT_upper}, stated below with some variables changed to avoid confusion with our notation.

\begin{Theorem}[Lemke Oliver--Thorne, {\cite[Theorem 1.2]{LOT_upper}}]
\label{thm:LOT_upper_2}
	Let $n \geq 2$.
	\begin{enumerate}
		 \item Let $a$ be the least integer for which $\binom{a+2}{2} \geq 2n + 1$. Then
		\[N_n(X) \ll X^{2a - \frac{a(a-1)(a+4)}{6n}}.\]
		\item Let $3 \leq b \leq n$ and let $a$ be such that $\binom{a+b-1}{b-1} > bn$. Then \[N_n(X) \ll X^{ab}.\]
	\end{enumerate}
\end{Theorem}

For a fixed superelliptic curve, i.e.\ choice of $m$ and $d$, we aim to find an integer $N$ such that for all $n \geq N$ satisfying $\gcd(m,d) \mid n$, the lower bound $N_{n,C}(X) \gg X^{\delta_n'}$ from Corollary \ref{cor:bound_improvements} holds. Below we summarize this procedure.
\begin{enumerate}
	\item Set $\varepsilon(n,m,d) = 1.564(\log n)^2$ and find $N_0$ such that $(**)$ is satisfied for all $n \geq N_0$. (Note that $(*)$ satisfied by Theorem \ref{thm:LOT_upper}.)
	\item Recalling $n_0 = \max(d, \Frob(m,d)+1)$, set
        \[n_0' = \begin{cases} n_0 & m \text{ is prime or } m \mid d,\\
            \max\left(n_0, \frac{m^2(m-1)}{2}  \right) & \text{otherwise}.
        \end{cases}\]
        Then for $n_0' \leq n \leq N_0$, use Theorem \ref{thm:LOT_upper_2} to search for $a$, $b$ values to find $\varepsilon(n,m,d)$ satisfying both $(*)$ and $(**)$.
\end{enumerate}
Note that $n \geq n_0'$ ensures that our parameterization strategy produces an irreducible polynomial family $F(t)$ by Proposition \ref{prop:param_irred_Sn}. For several small values of $m$ and $d$, we compute such an $N$ with this procedure, displayed below in Figure \ref{fig:improved_bound}.

\begin{center}
\begin{figure}[ht]
	\footnotesize
	\begin{tabular}{c || c | c || c | c || c | c || c | c || c | c || c | c || c | c }
		$m$ & \multicolumn{2}{c||}{2} & \multicolumn{2}{c||}{3} & \multicolumn{2}{c||}{4} & \multicolumn{2}{c||}{5} & \multicolumn{2}{c||}{6} & \multicolumn{2}{c||}{7} & \multicolumn{2}{c}{10} \\ \hline 
		$d$ & $n_0'$ & $N$ & $n_0'$ & $N$ & $n_0'$ & $N$ & $n_0'$ & $N$ & $n_0'$ & $N$ & $n_0'$ & $N$ & $n_0'$ & $N$  \\ \hline \hline
		3 & 3 & 106 & 3 & 552 & \multicolumn{2}{c ||}{} & \multicolumn{2}{c ||}{} & \multicolumn{2}{c ||}{} & \multicolumn{2}{c ||}{} & \multicolumn{2}{c}{}\\ \hline
		4 & 4 & 108 & 6 & 561 & 4 & 1164 & \multicolumn{2}{c ||}{} & \multicolumn{2}{c ||}{} & \multicolumn{2}{c ||}{} & \multicolumn{2}{c}{} \\ \hline
		5 & 5 & 110 & 8 & 563 & 24 & 1168 & 5 & 2015 & \multicolumn{2}{c ||}{} & \multicolumn{2}{c ||}{} & \multicolumn{2}{c}{} \\ \hline
		6 & 6 & 112 & 6 & 558 & 24 & 1162 & 20 & 2030 & 6 & 3192 & \multicolumn{2}{c ||}{} & \multicolumn{2}{c}{}\\ \hline
		7 & 7 & 114 & 12 & 573 & 24 & 1174 & 24 & 2034 & 180 & 3210 & 7 & 4438 & \multicolumn{2}{c}{}\\ \hline
		10 & 10 & 120 & 18 & 585 & 24 & 1166 & 10 & 2020 & 180 & 3196 & 54 & 4485 & 10 & 10860\\ \hline
		100 & 100 & 236 & 198 & 750 & 100 & 1256 & 100 & 2110 & 195 & 3376 & 594 & 5284 & 100 & 10940\\ \hline
		1000 & 1000 & 1000 & 1998 & 1998 & 1000 & 2040 & 1000 & 3045 & 1995 & 5074 & 5994 & 8892 & 1000 & 11800
	\end{tabular}

\caption{When is Corollary \ref{cor:bound_improvements} taking effect?}
\label{fig:improved_bound}
\end{figure}
\end{center}

\begin{Example}\label{ex:running_example_X_final}
    Let us return to our running concrete example (see Examples \ref{ex:running_example_I_S6point}, \ref{ex:running_example_II_S14point}, \ref{ex:running_example_III_higherdeg}, \ref{ex:running_example_IV_irr}, \ref{ex:running_example_V_galois}, \ref{ex:running_example_VI_countsetup_6}, \ref{ex:running_example_VII_countsetup_14}, \ref{ex:running_example_VIII_multiplicity}, \ref{ex:running_example_IX_result}),
    \[C \colon y^3 = x^4 + 1.\]        
    The improved bound of Corollary \ref{cor:bound_improvements} takes effect whenever $n \geq 561$. For $n=561$, we have
    \begin{align*}
        \delta_{561} &= 9668399/88121880 \approx 0.1097 \\
        \delta_{561}' &= 41851/376992 \approx 0.1110.
    \end{align*}
    This improvement puts $\delta_{561}'$ about 93\% of the way from $\delta_{561}$ to $1/9$, the limiting value of $\delta_n, \delta_n'$ as $n \to \infty$.
\end{Example}

\section{Geometric sources of higher degree points}\label{sec:failsafe}

Let $C$ be a superelliptic curve over $\Q$ given by an affine equation of the form $y^m=f(x)$ where $f(x)$ has degree $d$. The parametrization strategy in \eqref{eq:param} produces points on superelliptic curves that generate degree $n$ field extensions. The strategy fails to produce degree $n$ extensions when $\text{gcd}(m,d) \nmid n$ in general. In this section, we attempt to provide some heuristics for why one should expect degree $n$ points on superelliptic curves with $\text{gcd}(m,d) \nmid n$ to appear less often compared to degree $n$ points with $\text{gcd}(m,d) \mid n$. 

In the case of hyperelliptic curves, $m = 2$ and $\gcd(2, d) = 2$, this parametrization does not produce any odd degree points (cf. \cite{keyes}). This is consistent with a result of Bhargava--Gross--Wang \cite{BGW} which finds that for any genus $g \geq 2$, a positive proportion of everywhere locally soluble hyperelliptic curves have no odd degree points (and thus that a positive proportion of all genus $g$ hyperelliptic curves have no odd degree point). 

 While we are far from proving an analogous result to \cite{BGW} for degree $n$ points with $\text{gcd}(m,d) \nmid n$ on superelliptic curves, we attempt to give some heuristics and examples suggesting that points of degree $n$ with $\text{gcd}(m,d) \mid n$ appear more often than those with $\text{gcd}(m,d) \nmid n$ and we ask the following:

 \begin{question} What, if anything, can be said about the sparcity or abundance of various degrees $n$ of points on superelliptic curves given by \eqref{eq:sec}? In particular, can something be said in terms of the relationship of $n$ to the quantities $m$, $d$, and $\text{gcd}(m,d)$?
 \end{question}
 
Another way to phrase this question is in terms of the index of the curve $C/K$. The \textit{index} of a curve $C$, denoted $I(C)$, is the greatest common divisor of degrees $[L : K]$, where $L/K$ ranges over algebraic extensions such that $C(L)\neq \emptyset$. See \cite{GOLQL, Sharif} for more on the index of a curve. The result of Bhargava--Gross--Wang \cite{BGW} can be phrased as stating that a positive proportion of hyperelliptic curves over $\Q$ have index 2 over $\Q$. 

For a general superelliptic curve $C/\Q$, one can ask whether its index over $\Q$ is related to $\text{gcd}(m,d)$. It is already clear for instance that $I(C) \mid \text{gcd}(m,d)$ but we ask if more is true. If the exponent $m$ is prime, Creutz \cite{creutz} describes how descent can be used to determine that $\Pic^1(C)(\Q) = \emptyset$ --- which implies the index of $C/\Q$ is $m$ --- and gives a specific example of a curve with $m=3$, $d=6$ with index 3 (see \cite[Example 7.3]{creutz}). At present, the authors are not aware of similar explicit examples for other $(m,d)$ pairs or of families of superelliptic curves with index $\gcd(m,d)$ aside from $m=2$.

\subsection{Arithmetic from geometry}\label{sec:arith_from_geo}
A geometric source from which we can expect to find infinitely many points on $C$ are maps to $\mathbb{P}^1$. The most apparent of these are the natural maps of degree $m$ and $d$ from our curve $C$ to $\mathbb{P}^1$. That is, we can get infinitely many points by pulling back along the degree $m$ and degree $d$ maps to $\mathbb{P}^1$. Thus we know there are infinitely many degree $n$ points that are either multiples of $d$ or multiples of $m$. For other discussions on sources of infinitely many points on different types of curves, or more general curves, see \cite{AbramovichHarris, BourdonVirayEtAl, DebarreFahlaoui, HarrisSilverman, KadetsVogt, SV}. In another direction, for discussions on finiteness of points in certain degrees, see \cite{GM, Levin, SongTucker, Vojta}. 

In what follows, for $n$ the degree of the points and $g$ the genus of the curve, we discuss maps from $C$ to $\mathbb{P}^1$ in the case $n<g$ and in the case $n\geq 2g$.
 
\subsubsection{The case of $n < g$}

We first wish to characterize potential sources of infinitely many points on $C$ of degree $n < g$. Suppose further that the exponent $m$ is prime (we remark about the composite case below).

Define the $n$-th symmetric product of $C$ as usual by $\text{Sym}^n(C):= C^n/S_n$. The points of $\text{Sym}^n(C)$ correspond to effective degree $n$ divisors on $C$. We have a natural map
\[ \alpha: \text{Sym}^n(C) \to \Pic^n (C),\]
 defined by taking $D \mapsto [D]$.  $\Pic^n(C)$ is a $g$-dimensional variety (it is a torsor of the Jacobian of $C$), and the image $\alpha(\text{Sym}^n(C))$, often denoted by $W_n$, is a proper closed subvariety of $\Pic^n(C)$.

Suppose there exists a degree $n$ divisor class $[D_0]$, defined over $\Q$. Then $\Pic^n(C)$ is isomorphic to the Jacobian of $C$, denoted $J_C$, by the map $[D] \mapsto [D] - [D_0]$, and we extend the map $\alpha$ above to $J_C$ by composition with the isomorphism. In the case where $m$ is prime, by a result of Zarhin \cite[Theorem 1.2]{Zarhin} we have that for a generic $C$, $J_C$ is geometrically simple. That is, generically $J_C$ does not contain a translated proper abelian subvariety and therefore $\alpha(\text{Sym}^n(C))$ does not contain an abelian subvariety. 

By a theorem of Faltings \cite{faltings94}, this implies there are only finitely many points of $\alpha(\text{Sym}^n(C))$ and therefore only finitely many points of $\text{Sym}^n(C)$ that do not come from a $g^r_n$ on $C$. 

\begin{Theorem}[Faltings, \cite{faltings94}]Let $X$ be a closed subvariety of an abelian variety $A$, with both
defined over a number field $K$. Then the set $X(K)$ equals a finite union $\cup B_i(K)$, where each $B_i$
is a translated abelian subvariety of $A$ contained in $X$.
\end{Theorem}

In other words, generically, there are only finitely many points of $\text{Sym}^n(C)$ apart from those coming from the positive dimensional fibers of $\alpha$. We know that for some $n$ (namely, $n=m$ or $n$ a multiple of $m$) the map $\alpha$ must have positive dimensional fibers, because in particular the points of $\Sym^n (C)$ that are the result of pulling back points from maps from $C$ to $\P^1$ (e.g.\ a $g^1_m$) map to a point of $J_C$. This is because the Jacobian of $\P^1$ is trivial. However, the lack of a complete characterization of the positive dimensional fibers prevents us from concluding anything about finiteness of $C(K)$ in  certain degrees.

For hyperelliptic curves, there is a complete characterization of the positive dimensional fibers (see e.g., Arbarello--Cornalba--Griffiths--Harris, \cite{ACGH} page 13). Any effective degree $n$ divisor $D$ having positive rank on a hyperelliptic curve $H$ must contain a sub-divisor of the form $P + \iota(P)$ where $P$ is some point on $H$ and $\iota$ is the hyperelliptic involution. In other words, the only positive dimensional fibers of the map $\alpha$ when $C$ is a hyperelliptic curve are multiples of the $g^1_2$ (i.e. the only source of infinitely many points is pulling back along the degree $2$ map to $\mathbb{P}^1$).   Gunther--Morrow in \cite[Proposition 2.6]{GM} use this and argue as above to show that for  100\% of hyperelliptic curves $C$ (asymptotically as $g \to \infty$), $C$ has finitely many degree $n < g$ points that do not arise from pulling back a degree $n/2$ point of $\P^1$.

\begin{remark} In the case of $m$ composite, we no longer have that $J_C$ is geometrically simple, however work of Occhipinti--Ulmer \cite{OcchipintiUlmer} provides a useful understanding of the abelian subvarieties that appear in the Jacobian. More precisely, for a fixed polynomial $f(x)$ with $m=p_1^{a_1}\dots p_l^{a_l}$ (composite), the curve $C_m \colon y^m=f(x)$ has maps to other curves of the form $C_{m'} \colon y^{m/p_i^{b_i}}=f(x)$ where $1\leq b_i \leq a_i$ and $m':= m/p_i^{b_i}$. These maps between curves induce homomorphisms from the Jacobian $J_{C_m'}$ to $J_{C_m}$.  They define $J_m^{\text{new}}$ to be the quotient of $J_{C_m}$ by the sum of the images of these morphisms for all proper divisors $m'$ of $m$. $J_{C_m}$ is isogenous to the product of $J_{m'}^{\text{new}}$ with $m'$ ranging over all divisors of $m$. They show that for some sufficiently large $M$, $J_{M}^{\text{new}}$ does not contain any abelian subvarieties of dimension less than or equal to the genus of $C$.
\end{remark}

\subsubsection{The case of $n\geq 2g$}

For a fixed curve $y^m=f(x)$ where $f(x)$ has degree $d$, the parametrization in Proposition \ref{prop:param} produces infinitely many points of sufficiently large degrees $n$ divisible by $\gcd(m,d)$. Choose finitely many such points $P_1 \dots P_w$ of degrees $n_1 \dots n_w$ on $C$.

We now illustrate how one can use such points to produce a degree $n=\sum\limits_{i=1}^w n_i$ map to $\mathbb{P}^1$, that is, another source of infinitely many points of degree $n$. In this case $n$ will (by construction) be a multiple of $\text{gcd}(m,d)$.

To each point $P_i$, one can associate an element of $\text{Sym}^{n_i} (C)$ i.e., the effective degree $n_i$ divisors $D_i$ defined over $\mathbb{Q}$ corresponding to the Galois conjugates of $P_i$. Take $D:= D_1+\dots+ D_w$.
Let $w$ be a positive integer large enough such that $n \geq 2g$.
Using that $C$ is smooth and integral, we may identify Weil divisors with line bundles (see e.g., \cite{HartshorneAG}, II.6.16), and hence consider the line bundle $L(D)$, which is defined over $\mathbb{Q}$.
By Riemann--Roch (see e.g.,\ \cite{HartshorneAG}, IV.1.3), the line bundle $L(D)$ is basepoint free and has 
\[h^0(C,L(D)) = h^1(C, L(D)) + n + 1 - g \geq g + 1 \geq 2,\]
and so the sections of $L(D)$ define a map to $\mathbb{P}^1$. 
We may assume that the sections of $L(D)$ define a degree $n$ map to $\mathbb{P}^1$. If $h^0(C,L(D))$ is greater than $2$, we may instead take a sub-linear series.
Using a geometric version of the Hilbert Irreducibility Theorem (see e.g., \cite{Serreaspects} \S 9.2, Proposition 1), the fibers over all but a thin set of the rational points on $\mathbb{P}^1$ give us degree $n$ points on $C$. Note that if for our given curve, $\text{gcd}(m,d)=1$, then this produces a degree $n$ map to $\mathbb{P}^1$ giving us infinitely many points on $C$ for all $n$ sufficiently large.

\begin{remark} The above construction of a degree $n$ map to $\mathbb{P}^1$ began with points $P_1 \dots P_w$ coming from parametrization \eqref{eq:param} that each had degrees that were multiples of $\text{gcd}(m,d)$. The same construction could be carried out with $P_1 \dots P_{w+1}$ if one found a point $P_{w+1}$ on the curve not coming from the parametrization, but instead having some degree $n_{w+1}$ that is not a multiple of $\text{gcd}(m,d)$. The result of this would be that for $n$ sufficiently large, there is an infinite source of points that have degree $n$ (i.e. a degree $n$ map to $\mathbb{P}^1$) where $n$ is not a multiple of $\text{gcd}(m,n)$. 
\end{remark}

\begin{remark}
	If $g+1 \leq n < 2g$ and $L(D)$ is not basepoint free, we can still obtain a degree $n$ map to $\P^1$ from the curve minus the base point locus. By the “curve to projective” extension theorem, such a map extends to a map to $\P^1$ from the curve but the degree can be smaller by the degree of the base locus divisor. The degree of the base locus divisor must be divisible by the index of the curve.
\end{remark}

\subsection{Comparing points obtained via pullback to the parameterization}

A special case of the points produced by the parameterization, when $n=m$, our parametrization gives rise to points that come from pulling back along the $m$-to-one map to $\P^1$. The fields generated by these points generically have Galois group $C_m \rtimes C_{\varphi(m)}$, where $\varphi(m)$ is Euler's totient function.

To see this, consider for example, the curve $y^m=f(x)$, and set $\gamma(t)=\alpha_0$, $g(t)=t$, and $\eta(t)=h(t)=1$, this gives us the polynomial $F(t)=t^m-f(\alpha_0)$. We view this as a map from $\mathbb{P}^1$ to $\text{Sym}^m(C)$ sending the point $[\alpha_0:1]$ on $\P^1$ to the degree $n$ divisor consisting of the conjugates of the point $(\alpha_0, f(\alpha_0)^{1/m})$ on $C$. Thus the parametrization given in \eqref{eq:param} recovers degree $m$ points that come from pulling back along the $m$-to-one map to $\mathbb{P}^1$.

A description as above of the points coming from parametrization \eqref{eq:sec} as some copy of $\mathbb{P}^r$ (for some $r$) does not hold in general, so we ask the following question:
\begin{question} 
Is there a nice geometric characterization (perhaps as the rational points of a subvariety of $\Sym^n(C)$) of the points that arise from the parametrization given in \eqref{eq:param}?
\end{question}

\subsection{Heuristics for a special case using a result of Bhargava--Gross--Wang}

Suppose we have a curve $C$ given by an affine equation $y^m=f(x)$ where $f(x)$ has degree $d>4$. Suppose further that $m$ and $d$ satisfy  $2^i\mid \text{gcd}(m,d)$ where $i\geq 2$. Let $N=2k$ for $k$ an odd prime. In particular, for this case we have that $\gcd(m,d)\nmid N$. In what follows we find that many such curves $C$ have only finitely many degree $N$ points, subject to conditions on the size of $k$.

Let $C$ be the superelliptic curve given by $y^m=f(x)$ with $f(x)$ of degree $d$ and let $H$ be the hyperelliptic curve given by $y^2=f(x)$ (note that this is the same $f(x)$ as in the equation of $C$). We made the assumption that $d> 4$, so $H$ has genus at least $2$. We have a natural map $\phi$ from $C$ to $H$, given by sending points $\{ (x,\sqrt[m]{f(x)})\}$ to $\{ (x,\sqrt[2]{f(x)})\}$. If $P$ is a point of degree $N$ on $C$, we can map it to a point $P' = \phi(P)$ on $H$ as below. Let $\mathbb{Q}(P)$ and $\mathbb{Q}(P')$ be the extensions generated by a point $P$ on $C$ and by a point $P'$ on $H$, respectively.

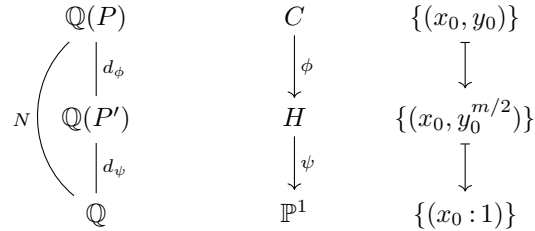
\[
\begin{tikzcd}
\mathbb{Q}(P) \arrow[no head]{d}{d_\phi} \arrow[dd, "N" ' , no head, bend right=49] & & C\arrow{d}{\phi} & \{ (x_0, y_0)\} \arrow[maps to]{d} \\
\mathbb{Q}(P') \arrow[no head]{d}{d_\psi} &  & H \arrow{d}{\psi} & \{(x_0,y_0^{m/2})\} \arrow[maps to]{d} \\
\mathbb{Q}                        &  & \mathbb{P}^1        & \{(x_0\hspace{.2em}\colon1)\}          
\end{tikzcd}
\]

By assumption, the degree of $[\mathbb{Q}(P) :\mathbb{Q}]$ is $N$. This means that the possibilities for $d_\psi$ and $d_\phi$ are as follows:\\
\begin{center}
\begin{tabular}{ |c||c|c|c|c| } 
\hline
Map & Case $1$ & Case $2$ & Case $3$ & Case $4$\\
\hline
\hline
$d_{\phi}$ & $1$ & $N$ & $k$ & $2$\\
\hline
$d_{\psi}$ & $N$ & $1$ & $2$ & $k$ \\
\hline
\end{tabular}
\end{center}
\vspace{3mm}

\textit{Case 1:} One should expect this to happen rarely as this would imply $\mathbb{Q}(P) = \mathbb{Q}(P')$, or equivalently $\Q\left(x_0, \sqrt[m]{f(x_0)}\right) = \Q\left(x_0, \sqrt{f(x_0)}\right)$, is an equality of degree $N$ number fields. \\

\textit{Case 2:} In this case $H$ has a rational point. Since we assumed $g(H)\geq 2$, Faltings' theorem \cite{faltings} implies that the set $H(\mathbb{Q})$ is finite. In fact, Shankar--Wang \cite{shankarwang} show that for even, monic hyperelliptic curves $H$ of genus $g(H) \geq 9$ with a marked rational non-Weierstrass point $\infty$, a positive proportion (tending to 100\% as $g(H) \to \infty$) have exactly two rational points, namely $\infty$ and $-\infty$, the conjugate of $\infty$ under the hyperelliptic involution. By assumption we have that $H$ is even, but even if $H$ is not monic, we may still be able to bound the number of rational points. Under certain technical assumptions (when $r\leq g(H)-3$, for $r$ the rank of $J_H$), Stoll \cite{Stoll} gives an explicit uniform bound for $\#H(\mathbb{Q})$ depending only on the genus of the curve and the rank of its Jacobian using the Chabauty--Coleman method \cite{Chabauty, Coleman} (see also \cite{McPoonen}). Therefore we may say that this does not happen often.\\

\textit{Case 3:} We have that $d_\phi$ is bounded above by the degree of $\phi$, so the Riemann--Hurwitz formula gives an upper bound 
\[d_{\phi} \leq \text{deg}(\phi) \leq \dfrac{g(C)-1}{g(H)-1}.\] Thus for $N$ sufficiently large (i.e.\ $k$ sufficiently large), Case 3 is excluded entirely.\\

\textit{Case 4:} Here $P'$ is an odd degree point of $H$. However, Bhargava--Gross--Wang \cite{BGW} show that a positive proportion of hyperellptic curves $H$ have no odd degree points, excluding this case. Note that for this positive proportion of curves, Case 2 also does not occur.\\

We conclude with an illustrative special case, in which we show that for many curves $C$ satisfying some conditions on $m, d, k$, we have at most finitely many points of degree $N$.

Let $f(x)$ be a squarefree polynomial of even degree $d = 2g + 2$. This gives a hyperelliptic curve with affine equation
\begin{equation}\label{eq:ours} y^2=f(x)=c_{2g+2}x^{2g+2}+ c_{2g + 1}x^{2g+1}+ \dots + c_{0} \end{equation}
with coefficients $c_i \in \Z$. We define the height of the polynomial $f(x)$ to be 
\[\height(f):= \text{max} \{ |c_i| \}.\]

We remark that Propositions $2.5$ and $2.6(2)$ of \cite{GM} hold for even degree hyperelliptic curves as in \eqref{eq:ours}. The results of Gunther--Morrow are stated for (odd) hyperelliptic curves with a rational Weierstrass point and their hyperelliptic curves are ordered by a slightly different height. We phrase the result in terms of densities of polynomials $f(x)$ so that our height is compatible with the height used in \cite{BGW}. We record the minor differences in the proofs in the following Lemma:

\begin{Lemma}\label{evenversion} 
Fix a genus $g \geq 2$ and an even integer $n < g$. Then for 100\% of squarefree polynomials $f(x)$ of degree $2g + 2$, ordered by height, the corresponding hyperelliptic curve $H/\Q$ given in \eqref{eq:ours} has at most finitely many degree $n$ points not obtained by pulling back degree $\frac{n}{2}$ points of $\mathbb{P}^1$.
\end{Lemma}

\begin{proof} 
First we show, as in \cite[Proposition 2.5]{GM}, that for $100\%$ of squarefree polynomials $f(x)$ the corresponding hyperelliptic curve $H$ has geometrically simple Jacobian. To see this, let $t_0, \dots t_{2g+2}$, so that we may show that polynomial $F(x, t_0, \dots, t_n)= t_{2g + 2}x^{2g+2}+\dots + t_{0}$ has Galois group $S_{2g+2}$ over $\mathbb{Q}(t_0, \dots t_{2g+2})$. Take the specialization with $t_2=\dots=t_{2g}=0$, $t_{0}=t_{1}=-1$, and $t_{2g+2}=1$. This gives us the polynomial $x^{2g+2}-x-1$, which is irreducible and has Galois group $S_{2g+2}$ by Corollary $3$ of \cite{Osada}.  By a result of Zarhin \cite{ZarhinHEC}, such $H$ have geometrically simple Jacobians. By the Hilbert Irreducibility Theorem \ref{thm:hilbert_irred} we see that $100\%$ of specializations of $F(x, t_1, \dots, t_n)$ have Galois group $S_{2g+2}$ and thus for $100\%$ of squarefree polynomials $f(x)$ the corresponding genus $g$ hyperelliptic curve $H$ with affine equation $y^2=f(x)$ has geometrically simple Jacobian.

The rest follows from exactly the same proof as \cite[Proposition 2.6]{GM}, outlined in Section \ref{sec:arith_from_geo}, except that since we do not assume there is a rational Weierstrass point, one defines an Abel Jacobi map $\Sym^n (H) \to J_H$ using a fixed degree $2$ divisor $D_0$ on the curve and sending $D\mapsto 2D-nD_0$ (we know such a divisor exists because of the map to $\mathbb{P}^1$).
\end{proof}

We also note that the above height on the polynomials $\height(f):= \max \{  |c_i| \}$ agrees with the height defined by Bhargava--Gross--Wang in \cite{BGW} when $H$ is embedded in the weighted projective space $\mathbb{P}(1,g+1,1)$ and expressed by the following equation:
\begin{equation}\label{eq:theirs} y^2=f(x,z)=c_{2g+2}x^{2g+2}+ c_{2g+1}x^{2g+1}z+ \dots + c_0z^{2g+2}. \end{equation}
Bhargava--Gross--Wang define the height of such a curve $H$ to be $\height'(H):= \text{max} \{ |c_i| \}$. This height on curves in weighted projective space corresponds exactly to the height $\height(f)$ on the defining polynomial $f(x)$ when we dehomogenize by taking $z=1$. Thus by \cite[Theorem 1]{BGW}, for a positive proportion of squarefree polynomials $f(x)$ ordered by height, $H$ has has no odd degree points.

\begin{Proposition}\label{prop:degN}
	Suppose $m$ and $d$ are positive even integers and $k$ is an odd prime satisfying
	\begin{itemize}
		\item $4 \mid m \mid d$,
		\item $\frac{m}{2} < k$,
		\item $N = 2k < \frac{d}{2} - 1$.
	\end{itemize}
	Then for a positive proportion of squarefree degree $d$ polynomials $f(x)$ ordered by height, the superelliptic curve given by $C \colon y^m = f(x)$ has finitely points of degree $N$. 
	
	Moreover, for such a curve $C$ and point $P \in C$ of degree $N$, the image $\phi(P) \in H$ as defined above is of degree $N$ and is not the pullback of a degree $k$ point on $\P^1$.
\end{Proposition}

\begin{proof}
	Let $H$ be the corresponding hyperelliptic curve with equation $H \colon y^2 = f(x)$. By methods of \cite{GM} (see Lemma \ref{evenversion}), we have that for $100 \%$ of polynomials $f(x)$ of degree $d$, the corresponding hyperelliptic curve has only finitely many points of degree $n < g(H)$ that are not the pullback of a degree $\frac{n}{2}$ point on $\mathbb{P}^1$. We also know that a positive proportion of such hyperelliptic curves do not have any odd degree points by \cite[Theorem 1]{BGW}. Thus for a positive proportion of polynomials $f(x)$, the hyperelliptic curve $H$ has both of these properties. For such $H$, let $C\colon y^m = f(x)$ be the superelliptic curve with map to $H$ given by $\phi \colon (x_0, y_0) \mapsto (x_0, y_0^{m/2})$ as above. Take $P$ to be a point on $C$ of degree $N$ with $\phi(P)=P'$ its image in $H$. By considering Cases 1 --- 4 above, we show there are only finitely many such $P$.
	
	Cases 2 and 4 are excluded by the fact that $H$ has no odd degree points. To see that Case 3 is impossible, we recall $k = d_\phi \leq \deg(\phi) = \frac{m}{2}$. This contradicts the hypothesis, so $d_\phi$ cannot be equal to $k$.
	
	All that remains is Case 1, in which both $P$ and its image $P'$ are degree $N$ points. Suppose $P = (x_0, y_0)$, so $P' = (x_0, y_0^{m/2})$, and the image of $P, P'$ in $\mathbb{P}^1$ is $x_0$. If $P'$ is the pullback of a degree $k=N/2$ point of $\mathbb{P}^1$, then $[\mathbb{Q}(x_0): \mathbb{Q}] = k$  and $f(x_0)$ is not a square in $\mathbb{Q}(x_0)$. However, this implies that the degree of $\sqrt[m]{f(x_0)}$ over $\mathbb{Q}(x_0)$ is greater than $2$, which contradicts that the degree of $P$ is $N$.
	
	Thus we see that $P'$ cannot be the pullback of a degree $k$ point of $\mathbb{P}^1$. Since $N < \frac{d}{2} - 1 = g(H)$ by our assumption, we have that only finitely many such $P'$ can exist. Hence at most finitely many points $P$ on $C$ of degree $N$ can exist.
\end{proof}

The argument for Case 1 in the proof of Proposition \ref{prop:degN} can be refined to prove Proposition \ref{prop:degN_2adic}, at the expense of the description of the image of $P$ in $H$. For this, we do not use \cite[Theorem 1]{BGW}, allowing us to obtain a proportion approaching 100\%.

\begin{Proposition}[See Proposition \ref{prop:degN_2adic}]\label{Prop:degN_2adic}
		Suppose $m, d$ are positive even integers such that $d > 4$. Let $N < \frac{d}{2} - 1$ have $2$-adic valuation strictly less than that of $m$, i.e.\ $v_2(N) < v_2(m)$. Then for a positive proportion approaching 100\% of squarefree degree $d$ polynomials $f(x)$, ordered by height, the superelliptic curve $C \colon y^m = f(x)$ has finitely many points of degree $N$.
\end{Proposition}

\begin{Example}
    Let $C$ be a superelliptic curve with $m=4$ and $d=24$. Generically, this will have genus $33$. The proposition above gives us that for a positive proportion of polynomials $f(x)$, there are finitely many points of degree $N=10,9,7,6,5,3,2,1$ on the curve $C\colon y^4=f(x)$. Note that, applying \cite[Corollary 0.3]{Vojta} to such $C$, we get finitely many points $P$ of degree $N \leq 9$ with the property that $P$ and $\psi(\phi(P))$ (i.e.\ the image under the usual degree $m$ map $C \to\P^1$) have the same field of definition.
    
   Specializing \cite[Corollary 0.3]{Vojta} in the case of superelliptic curves more generally gives information about a different but non-disjoint set of degrees as our Proposition \ref{Prop:degN_2adic}. Our result gives information about larger degrees $N$ when $d$ is large relative to $m$, but subject to the 2-adic valuation condition and only for a positive proportion of $f(x)$.
\end{Example}

\bibliographystyle{siam}
\bibliography{sec_bibliography}

\end{document}